\magnification 1200
\baselineskip 14pt
\input amssym.def
\input amssym.tex

%Numerazione delle proposizioni, diagrammi,...

\def\universale{(0.1)} %Succesione universale in G(1,3)

\def\canonicoincidenza{(1.1)} %Formula del canonico
\def\classeramificazione{(1.2)} %Formula della ramificazione
\def\equno{(1.3)} %Prima equazione per la focale
\def\eqdue{(1.4)}
\def\eqtre{(1.5)}
\def\eqquattro{(1.6)}
\def\classefocale{1.7} %Calcolo della classe della focale
\def\singfocale{1.8} %Nodale e cuspidale della focale
\def\fibratoSemple{(1.9)} %Diagramma che definisce il fibrato di Semple
\def\classeIX{(1.10)} %Classe di IX dentro il prodotto di P^3 per X
\def\classeDX{(1.11)}
\def\classeDDX{(1.12)}
\def\duedue{1.13} %Esempio (2,2) g=1
\def\duetre{1.14} %Esempio (2,3) g=1
\def\tretre{1.15} %Esempio (3,3) g=2

\def\invarianticorde{2.1}
\def\tangenteChow{2.2} %non lo so
\def\tangentecorde{2.3} %non lo so
\def\duesei{2.4} %Bisecanti a quartica ellittica

\def\univ{(3.1)} %Diagramma che definisce Q
\def\classi{3.2} %Classi delle bitangenti e dei flessi
\def\invarianti{3.3} %Invarianti delle bitangenti e dei flessi
\def\retta{3.4} %Retta contenuta in superficie focale
\def\bitangenti{3.5} %Invarianti della congruenza di bitangenti
\def\classiChern{3.6} %Classi di Chern dei prodotti simmetrici

\def\tangentecomplesso{4.1} %Spazio tangente al complesso tangente
\def\matrice{(4.2)} %Matrice di Y
\def\tangentebitangenti{4.3} %Tangente alla congruenza bitangenti
\def\tangentepolinomi{4.4} %Tang. ai polinomi con due radici doppie
\def\singbitangenti{4.5} %Singolarita' della congruenza bitangenti
\def\tangenteflessi{4.6} %Spazio tangente congruenza flessi
\def\singflessi{4.7} %Singolarita' della congruenza di flessi
\def\rigatabitangenti{4.8}
\def\rigataflessi{4.9}

\def\bigrado{5.1} %Bigrado della cong. di bitangenti a sup. singolare
\def\Veronese{5.2} %Esempio della Veronese
\def\seidue{5.3} %Esempio della duale della (2,6)
\def\dualWelters{5.4} %Esempio della duale della Welters
\def\Kummer{5.5} %Esempio della congruenza di Kummer
\def\focaleriducibile{5.6} %Congettura focale riducibile
\def\focalenonridotta{5.7} %Congettura focale ridotta
\def\bitangentiliscia {5.8} %Congettura cong. bitangenti liscie
\def\flessiliscia{5.9} %Congettura cong. flessi liscie
\def\bitangentifocale{5.10} %Congettura cong.bitangenti = congruenza

%Numerazione per la bibliografia

\def\ArrondoGross{[1]}
\def\Arkiv{[2]}
\def\Ciro{[3]}
\def\Fano{[4]}
\def\Goldstein{[5]}
\def\GH{[6]}
\def\Gross{[7]}
\def\Hudson{[8]}
\def\Johnsen{[9]}
\def\Schubert{[10]}

\def\McCrory{[12]}
\def\McCrorybis{[13]}
\def\Peskine{[14]}
\def\Rothuno{[15]}
\def\Rothdue{[16]}
\def\Salmon{[17]}
\def\Schumacher{[18]}
\def\Verra{[19]}
\def\Walker{[20]}
\def\Welters{[21]}

\def\qed{\hfill\vbox{\hrule\hbox{\vrule\kern3pt
    \vbox{\kern6pt}\kern3pt\vrule}\hrule}} %Fine di dimostrazione

\def\mapright#1{
     \smash{\mathop{\longrightarrow}\limits^{#1}}}

%Freccie con nome

 2
\font\large=cmbx10 scaled \magstep 2
 4

\font\semilarge=cmbx10 scaled \magstep 1

\def\qed{\hfill\vbox{\hrule\hbox{\vrule\kern3pt
    \vbox{\kern6pt}\kern3pt\vrule}\hrule}} %Fine di dimostrazione

\def\mapright#1{\smash{
   \mathop{\longrightarrow}\limits^{#1}}}

 %Freccie con nome

\centerline{\large A FOCUS ON FOCAL SURFACES} \centerline{E.
Arrondo, M. Bertolini and C. Turrini}
\bigskip
\bigskip

Many classical problems in algebraic geometry have regained interest when
techniques from differential geometry were introduced to study them. The
modern foundations for this approach has been given by Griffiths and
Harris in \GH, who obtained in this way several classical and new results
in algebraic geometry. More recently, this idea has been successfully
followed by McCrory, Shifrin and Varley in \McCrory\ and \McCrorybis\ to
study differential properties of hypersurfaces in ${\Bbb P}^3$ and ${\Bbb
P}^4$. In fact these two papers have greatly influenced the present work.

In this spirit, the subject of this paper is the systematic study of focal
surfaces of smooth congruences of lines in ${\Bbb P}^3$. This is indeed a
clear example of a topic of differential nature in algebraic geometry. The
study of such congruences has been very popular among classical
algebraic geometers one century ago. Especially Fano has given many 
important contributions to this field. An essential ingredient in his work
has been the focal surface of the congruence. This point of view has been
retaken by modern algebraic geometers, such as Verra and  Goldstein, and
also by Ciliberto and Sernesi in higher dimension.

What we find amazing in the papers by the classics is how much
information they were able to provide about the focal surface of
the known examples of congruences, in particular about its singular locus
(and more especially about fundamental points). They seemed to have in mind
some numerical relations that they never formulated explicitly. And even
nowadays such kind of relations would require deep modern techniques, like
multiple-point theory, but also this powerful machinery is not a priori
enough since some generality conditions need to be satisfied.

As a sample of this, the degree and class of the focal surface --the only
invariants easy to compute-- can be derived immediately from the
Riemann-Hurwitz formula. However these invariants, even in the 
easiest examples (see Example \duesei\ or Remarks after Corollary
\singflessi) seem to be wrong at a first glance. This is due to the
existence of extra components of the focal surface or to the possibility
that the focal surface counts with multiplicity, although this was never
mentioned explicitly by the classics. Even in \Goldstein, these
possibilities seem not to have been considered.

The starting point of this work was to understand how the classics
predicted the number of fundamental points of a congruence. We
only know of one formula in the literature involving this number, which is
however wrong (see Example \tretre\ and the remark afterwards). So our
first goal was to use modern techniques in order to rigorously obtain
some of the classical results on the topic. Specifically, by regarding the
focal surface of a congruence as a scheme, we reobtain its invariants
(degree, class, class of its hyperplane section, sectional genus, and
degrees of the nodal and cuspidal curves) and give them a precise sense. 

We also restrict our attention to congruences of bisecants to a curve, or
flexes to a surface (since they are special cases in the work by
Goldstein), or bitangents to a surface (since all the lines of a
congruence are bitangent to the focal surface). In particular, we prove
that no congruence of flexes to a smooth surface is smooth, and that a
congruence of bitangents to a smooth surface is smooth if and only if the
surface is a quartic not containing any line. Another important reason to
study these types of congruences is that their focal surfaces have the 
unexpected or multiple components mentioned above. We give a precise
geometrical description of these components and also conjecture that these
congruences are the only ones for which the focal surface has such a
behavior.

In order to obtain all the above results, we combine a local differential
analysis with global methods from intersection theory. In fact,
we consider that many of the techniques we develop are interesting by
themselves.

In section \S0, we give the basic definitions about congruences and
their focal surfaces. In section \S1, we obtain the classical
invariants of the focal surface. The key new technique in this
section is to use the construction given in \Arkiv\ of varieties
parametrizing infinitely close points of a given variety.

In sections \S2, \S3 and \S4 we obtain all the invariants of the
congruences given by bisecants to a smooth curve or bitangents or
flexes to a smooth surface in ${\Bbb P}^3$. In these sections we
again adapt some natural constructions to our setting. For instance
the constructions at the beginning of section \S3 are clearly
influenced by the ones in \McCrory.

Finally, section \S5 is devoted to relate the behavior of
congruences of bitangents to a smooth surface to the behavior of
general congruences. We give there several examples and conjectures
of what we expect to happen in general.

\bigskip

\noindent ACKNOWLEDGMENTS: Most of this work has been made in the
framework of the Spanish-Italian project HI1997-123. Partial
support for the first author has been provided by DGICYT grant
PB96-0659, while the two last authors were supported by the Italian
national research project ``Geometria Algebrica, Algebra
Commutativa e Aspetti Computazionali'' (MURST cofin. 1997). We want
also to acknowledge the extremely useful help that has been for us
the extensive use we did of the Maple package Schubert (\Schubert).
We also had the invaluable help of Mar\'{\i}a Jes\'us
V\'azquez-Gallo, who kindly adapted to our setup the sophisticated
Maple package she created for making computations in the Chow ring
of several parameter spaces. We finally thank Trygve Johnsen for
kindly giving us the reference for the formula about stationary
bisecants we were looking for.

\bigskip
\bigskip

\noindent{\semilarge \S0. Notations and definitions.}
\bigskip

We will work over an algebrically closed field of characteristic zero.
We will denote by $G(1,3)$ the Grassmann variety of lines in ${\Bbb
P}^3$. If $I\subset {\Bbb P}^3\times G(1,3)$ is the incidence
variety of pairs $(x,L)$ such that $x$ is a point of the line $L$,
then any of the projections $p_1$ and $p_2$ provides $I$ with a
structure of projective bundle. In fact, $I={\Bbb P}(\Omega_{{\Bbb
P}^3}(2))$ (where ${\Bbb P}$ will always mean for us the space of
rank-one quotients), and the tautological quotient line
bundle is just the pull-back of the hyperplane line bundle on
$G(1,3)$ (considered as a smooth quadric in ${\Bbb P}^5$). On the
other hand, if we consider the Euler sequence on ${\Bbb P}^3$
$$0\to \Omega_{{\Bbb P}^3}(1)\to
H^0({\Bbb P}^3,{\cal O}_{{\Bbb P}^3}(1))\otimes{\cal O}_{{\Bbb P}^3}
\to{\cal O}_{{\Bbb P}^3}(1)\to 0$$
\noindent and pull it back to $I$ via $p_1$ and then push it down to
$G(1,3)$ via $p_2$ we get the universal exact sequence on
$G(1,3)$
$$\leqno{\universale}\ \ \ \ \ \ \ \ \ \ \ \ \ \ \ \ \ \ \ \ \
0\to S^*\to H^0({\Bbb P}^3,{\cal O}_{{\Bbb P}^3}(1))\otimes{\cal
O}_{G(1,3)}\to Q\to 0.$$
\noindent Here $S$ and $Q$ are the rank-two universal vector
bundles, and $I$ can also be viewed as ${\Bbb P}(Q)$.

Given a point $x\in{\Bbb P}^3$, we define the {\it alpha-plane}
associated to it as the set $\alpha(x)\subset G(1,3)$ of all lines in
${\Bbb P}^3$ passing through it. Similarly, given a plane
$\Pi\subset {\Bbb P}^3$, we define the {\it beta-plane} associated
to it as the set $\beta(\Pi)$ of all lines in ${\Bbb P}^3$ contained
in $\Pi$. If $x\in\Pi$, we will write $\Omega(x,\Pi)$ for the pencil
of lines contained in the plane $\Pi$ and passing through the point
$x$.

By {\it congruence} we will mean a surface $X\subset G(1,3)$. Any
congruence $X$ has a {\it bidegree} $(a,b)$, where $a$ (called the
{\it order of the congruence}) is the intersection number of $X$ with
an alpha-plane, and $b$ (called the {\it class of the congruence}) is
the intersection number with a beta-plane. Equivalently,
$a=c_2(Q_{|X})$, and $b=c_2(S_{|X})=c_1(Q_{|X})^2-c_2(Q_{|X})$.

A congruence can be regarded (under the Pl\"ucker embedding of
$G(1,3)$) as a surface contained in a smooth quadric of ${\Bbb P}^5$.
In particular, we can define the {\it sectional genus} of a
congruence as the genus of the curve obtained by intersecting the
surface with a hyperplane of ${\Bbb P}^5$. We will usually denote it
with $g$.

A line in ${\Bbb P}^3$ can also be viewed as a line in the dual ${{\Bbb
P}^3}^*$, so that a congruence $X\subset G(1,3)$ induces another
congruence $X^*\subset G(1,{{\Bbb P}^3}^*)$, which we will call the
{\it dual congruence} of $X$. It is clear that, if $X$ has degree
$(a,b)$ then $X^*$ has bidegree $(b,a)$. A congruence and its dual
have the same sectional genus (in fact both Pl\"ucker embeddings are
naturally isomorphic).

If we restrict the above projections $p_1$ and $p_2$ to
$I_X:=p_2^{-1}(X)$ then we get a map $q_X: I_X\to{\Bbb P}^3$ which
is generically $a:1$ and a map $p_X: I_X\to X$. We have the
following definitions:

\noindent{\bf Definitions:} A point $x\in{\Bbb P}^3$ is a {\it
fundamental point} of $X$ if $q_X^{-1}(x)$ is not a finite set.
Dually, a {\it fundamental plane} of $X$ is a fundamental point of
$X^*$, i.e. a plane containing infinitely many lines of the congruence.
The {\it focal locus} of $X$ is the branch locus (typically a surface)
of $q_X$. The elements of the focal locus are called {\it focal points}
of $X$. Dually, a {\it focal plane} of $X$ is a focal point of $X^*$.
Equivalently (\Goldstein\ Lemma 4.4), a focal point $x\in{\Bbb P}^3$
is characterized by the fact that there exists a line $L$ of the
congruence such that the embedded tangent plane of $X$ at the point
represented by $L$ meets the alpha plane $\alpha(x)$ in at least a
line of ${\Bbb P}^5$ (i.e. a pencil of lines of ${\Bbb P}^3$). This is
in fact the definition of focal point given by Goldstein.

If we write $H$ and $K$ respectively for the classes of the hyperplane
section and the canonical divisor of $X$, and $h$
for the class of the hyperplane section of ${\Bbb P}^3$, it is not
difficult to see that $c_1(T_{I_X})=2h-K-H$, so that the class of the
ramification locus of $q_X$ is $2h+K+H$. In particular, we obtain from
here the very well-known result that a general line $L$ of a
congruence contains two focal points $x_1,x_2$ (counted with
multiplicity) such that $(x_1,L)$ and $(x_2,L)$ lie in the
ramification locus of $q_X$.

\noindent{\bf Definition:} We will call a {\it focal line} $L$ of a
congruence to a line of a congruence such that all of its points are
focal. Again from \Goldstein\ Lemma 4.4, this means that the
embedded tangent plane of $X$ at the point represented by $L$ meets in
a pencil all the alpha planes $\alpha(x)$ for which $x\in L$. Then, a
line $L$ is focal if and only if its embedded tangent plane is a
beta-plane.

Let $X_0$ be the open set of non-focal lines of a congruence $X$.
Then the restriction of the map $p_X^{-1}(X_0)\to X_0$ to the
ramification locus of $q_X$ is finite (typically of degree two, but
it could happen a priori that any line contains only one focal
point counted twice). Hence, the branch locus of this restriction
has at most two components.

\noindent {\bf Definition:} We will call the {\it strict focal
surface} of a congruence $X$ to the closure $F_0$ of the reduced
structure of the branch locus of $p_X^{-1}(X_0)\to {\Bbb P}^3$. To
distinguish from this, we usually refer to the focal locus $F$ (as a
scheme) as the {\it total focal surface}.

\noindent{\bf Remarks:} 1) We abused the notation in the above
definition. First of all, $X_0$ could be empty. As observed in the
definition of focal line, this would imply that the embedded
tangent plane of $X$ at any point is a beta-plane. In this case,
the congruence itself is a beta-plane (see for instance \Goldstein,
Corollary 4.5.1). On the other hand, $F_0$ could be either a point
(and then $X$ is an alpha-plane) or a curve, which would mean that
$X$ is the congruence of bisecants to that curve. As we will
observe later, such a congruence is only smooth when the curve is a
twisted cubic or an elliptic quartic.

2) It is not superfluous to take the reduced structure in the above
definition. As we will see in section $3$, the focal locus can
appear with high
multiplicity for congruences of bitangents or flexes to a surface in
${\bf P}^3$.

3) The total focal surface could have more components different
from $F_0$ when $X\setminus X_0$ is a curve. Such a curve will have
the property that its embedded tangent line at each point is
contained in $G(1,3)$, so that the corresponding extra components
of the focal surface will be developable ruled surfaces or cones.
The existence of these ruled surfaces seems not to have considered
by Goldstein. In fact, the number of components of the focal
surface can be bigger than two, as will be shown in Example
\seidue. \bigskip

We end this section of background definitions and results by recalling
a classical invariant for surfaces in ${\Bbb P}^3$ that we will use
frequently:

\noindent{\bf Definition:} If $\Sigma\subset{\Bbb P}^3$ is a
surface, we will write $\mu_1$ for the class of its hyperplane
section. It is clear that a surface and its dual have the same
invariant $\mu_1$.

\bigskip
\bigskip

\noindent {\semilarge \S 1. Numerical invariants of the focal surface
of a smooth congruence.}
\bigskip

Along this section, $X\subset G(1,3)$ will be a smooth congruence of
lines in ${\Bbb P}^3$, $H$ and $K$ will be the hyperplane and
canonical classes respectively, and $F$ will be the total focal
surface of $X$.

In order to better understand the geometry and the numerical invariants
of $F$ (in particular $\mu_1$), it is convenient to work in the
complete flag variety of points, lines and planes rather than only in
the incidence variety of points and lines. We consider then
$$A_X:=\{(x,L,\Pi)\in {\Bbb P}^3\times X\times {{\Bbb P}^3}^*\ | \
x\in L\subset\Pi\}.$$
Let $q_{13}:A_X\to J\subset{\Bbb P}^3\times{{\Bbb P}^3}^*$ and
$q_2:A_X\to X$ be the obvious projections, $J$ being the incidence
variety of points and planes. Our goal is to directly obtain the
focal variety in $J$, so that we construct simultaneously its dual.
For this purpose, we analyze the ramification locus of $q_{13}$.
First, we observe that the map $q_2$ factors $A_X\to I_X\mapright{p_X}
X$. The second morphism is the restriction of the projective bundle
$p_2:I\to G(1,3)$, so that $I_X\cong{\Bbb P}(Q_{|X})$, and the tautological line
bundle is the pullback of the hyperplane section $h$ of ${\Bbb P}^3$.
Similarly, the first morphism is a projective bundle and $A_X\cong
{\Bbb P}(p_X^*S_{|X})$, and its tautological line bundle is the
pullback of the hyperplane section $h^*$ of ${{\Bbb P}^3}^*$. From
this, it is not difficult to compute the Chern classes of $T_{A_X}$:
$$\leqno{\canonicoincidenza}\ \ \ \ \ \ \ \ \ \ \ \ \ \ \ \ \ \ \ \ \
\ \ \ \ \ \ c_1(T_{A_X})=2h+2h^*-2H-K$$
$$c_2(T_{A_X})=h^2+4hh^*+{h^*}^2-3hH-3h^*H
-2hK-2h^*K+2H^2+2HK+c_2(T_X)$$
On the other hand, the Chern classes of the incidence variety $J$
are:
$$c_1(T_J)=3h+3h^*$$
$$c_2(T_J)=3h^2+10hh^*+3{h^*}^2$$
Hence the class in $A_X$ of the ramification locus $R$ of $q_{13}$
will be, using Porteous formula ($R$ will be a surface, since $A_X$
has dimension four, and $J$ has dimension five),
$$\leqno{\classeramificazione}\ \ \ \ \ \ \ \
[R]=2hh^*+hH+h^*H+hK+h^*K+2H^2+2HK+K^2-c_2(T_X).$$ Then, for a
general line $L$ in the congruence, one expects to find two elements
$(x_1,L,\Pi_1)$, $(x_2,L,\Pi_2)$ in $R$, and it seems reasonable to
think that $x_1$, $x_2$ are focal points for $X$, that
$\Pi_1$, $\Pi_2$ are focal planes and that each $\Pi_i$ is the tangent
plane of the focal surface at $x_i$ (it is a very well-known result
that the set of focal planes is the dual of the focal surface).
However, the last of the statements is not true, but $\Pi_1$ is the
tangent plane of the focal surface at $x_2$ and reciprocally $\Pi_2$
is the tangent plane at $x_1$. Let us check this in local coordinates.

Fix an element $(x,L,\Pi)$ in $R$ and choose coordinates
$z_0,z_1,z_2,z_3$ in ${\Bbb P}^3$ so that $x$ is the point of
coordinates $(1:0:0:0)$, $L$ is the line $z_2=z_3=0$ and $\Pi$ is the
plane $z_3=0$. We can take $u,v$ to be a system of parameters of $X$
at $L$ and assume that near $L$ the lines of the congruence are given
by the span of the rows of the matrix
$$\pmatrix{
1&0&f&g\cr 0&1&h&k}$$ where $f,g,h,k$ are regular functions in a
neighborhood of $L$. We can take then a system of coordinates
$\lambda,u,v,\mu$ for $A_X$ near $(x,L,\Pi)$ to represent the point
$x(\lambda,u,v)=(1:\lambda:f+\lambda h:g+\lambda k)$ inside the
above line $L(u,v)$ and the plane $\Pi(u,v,\mu)$ containing them of
equation $z_3+\mu z_2=(g+\mu f)z_0+(k+\mu h)z_1$. On the other
hand, we can take affine coordinates $a_1,a_2,a_3$ to represent the
points $(1:a_1:a_2:a_3)\in{\Bbb P}^3$ and affine coordinates
$u_0,u_1,u_2$ to represent the plane $z_3-u_2z_2=u_0z_0+u_1z_1$. We
could remove one coordinate to work in $J$, locally defined as
$a_3-u_2a_2=u_0+u_1a_1$, but we prefer to keep the symmetry.
Therefore a local expression for $q_{13}$ is given by
$$(\lambda,u,v,\mu)\mapsto
(\lambda,f+\lambda h,g+\lambda k,g+\mu f,k+\mu h,\mu).$$
Its Jacobian matrix is then
$$\pmatrix{
1&       h       &       k       &     0     &     0     &0\cr
0&f_u+\lambda h_u&g_u+\lambda k_u&g_u+\mu f_u&k_u+\mu h_u&0\cr
0&f_v+\lambda h_v&g_v+\lambda k_v&g_v+\mu f_v&k_v+\mu h_v&0\cr
0&       0       &       0       &     f     &     u     &1}$$
We immediately see that this matrix has not maximal rank if and only
if the two middle rows are linearly dependent. Since the four columns
of this submatrix are linearly dependent, the local equations of $R$
are:
$$\leqno{\equno}\ \ \ \ \ \ \ \ \ \ \ \ \ \ \ \ \ \ \ \ \ \ \ \ \ \ \
\left|\matrix{
f_u+\lambda h_u&g_u+\lambda k_u\cr
f_v+\lambda h_v&g_v+\lambda k_v}\right|=0$$
$$\leqno{\eqdue}\ \ \ \ \ \ \ \ \ \ \ \ \ \ \ \ \ \ \ \ \ \ \ \ \ \ \
\left|\matrix{
g_u+\mu f_u&k_u+\mu h_u\cr
g_v+\mu f_v&k_v+\mu h_v}\right|=0$$
$$\leqno{\eqtre}\ \ \ \ \ \ \ \ \ \ \ \ \ \ \ \ \ \ \ \ \ \ \ \ \ \ \
\left|\matrix{
f_u+\lambda h_u&k_u+\mu h_u\cr
f_v+\lambda h_v&k_v+\mu h_v}\right|=0$$
Equation \equno\ means that the value of $\lambda$ is so that
$x(\lambda,u,v)$ is a focal point in $L(u,v)$, while \eqdue\ means that
$\Pi(u,v,\mu)$ is a focal plane. For a ``general'' value of $u,v$ there
would be two possible values of $\lambda$ and
$\mu$, and \eqtre\ should be interpreted as a way of
assigning to each of the two focal points in the line one of the two
focal planes. The key observation is that, substracting \equno\
multiplied by $\left|{h_u\atop h_v} {k_u\atop k_v}\right|$ and
\eqdue\ multiplied by $\left|{f_u\atop f_v}{h_u\atop h_v}\right|$
one gets
\eqtre\ multiplied by:
$$\leqno{\eqquattro}\ \ \ \ \ \ \ \ \ \ \ \ \ \ \ \ \ \ \ \ \ \ \ \ \ \
\
\left|\matrix{
f_u+\lambda h_u&k_u-\mu h_u\cr
f_v+\lambda h_v&k_v-\mu h_v}\right|=0.$$
This means that \eqquattro\ is the other way of assigning to each focal
point a focal plane (and we want to prove that this is the ``right''
one).

Assume now for simplicity that we chose $L$ containing exactly two focal
points $x,x'$ and contained in two focal planes $\Pi,\Pi'$. Then there
are two corresponding local expressions $\lambda,\lambda'$ in terms of
$u,v$ verifying \equno\ and two local expressions $\mu,\mu'$ verifying
\eqdue, and such that each of the pairs $(\lambda,\mu)$ and
$(\lambda',\mu')$ verify \eqtre, while the pairs $(\lambda,\mu')$ and
$(\lambda',\mu)$ verify \eqquattro. In particular, the assignement
$$(u,v)\mapsto (\lambda,f+\lambda h,g+\lambda k)$$
is a local parametrization of the focal surface near $x$. However, the
tangent plane at it is not $\Pi$, but $\Pi'$. Indeed, let $\mu'_0$ the
nonzero solution of \eqdue\  for $u=v=0$. Then $\Pi'=\Pi(0,0,\mu'_0)$
has equation $z_3+\mu'_0 z_2=0$. To check that $\Pi'$ is tangent we
need to show that, substituting the above parametrization in the
equation of $\Pi'$ we do not get linear terms. The substitution becomes
$z_3+\mu'_0 z_2=g+\lambda k+\mu'_0(f+\lambda h)$, and we need to check
that the partial derivatives vanish at $u=v=0$ (and hence also
$\lambda=0$). These partial derivatives are:
$$g_u(0,0)+\mu'_0f_u(0,0)\ {\rm and}\ g_v(0,0)+\mu'_0f_v(0,0)$$
To check this vanishing we first observe that \equno\ for $u=v=0$
implies that each vanishing implies the other. On the other hand,
\eqquattro\ implies that
$$\left|\matrix{f_u(0,0)&h_u(0,0)\cr f_v(0,0)&h_v(0,0)}\right|\mu'_0=
\left|\matrix{h_u(0,0)&g_u(0,0)\cr h_v(0,0)&g_v(0,0)}\right|.$$
{}From this it is easy to conclude that $\Pi'$ is indeed the tangent
plane.

We can now use the above calculations to prove the following:

\proclaim Proposition \classefocale. Let $X$ be a smooth
congruence, let $F$ be its total focal surface and consider $\tilde
F\subset A_X$ to be the closure of the set of elements $(x,L,\Pi)$
such that $(x,L)$ is a ramification point of $q_X$ and $\Pi$ is the
tangent plane to $F$ at a smooth point $x$. Then:
\item{1)} $\tilde F$ is linked to $R$ in the complete intersection of
the pullbacks to $A_X$ of the ramification loci of $q_X$ and
$q_{X^*}$. In particular, the cycle class of $\tilde F$ in $A_X$ is
$$[\tilde F]=2hh^*+hH+h^*H+hK+h^*K-H^2+c_2(T_X).$$
\item{2)} The focal surface $F$ has degree $2a+2g-2$ and, if it is
reduced, has class $2b+2g-2$, $\mu_1=a+b+4g-4-K_X^2+12\chi({\cal
O}_X)$, sectional (geometric) genus $9g-8-b+K_X^2$ and $\chi({\cal
O}_{\tilde F})=6g-6-a-b+K_X^2+2\chi({\cal O}_X)$.

{\it Proof:} The fact that $\tilde F$ and $R$ are linked is just
the geometrical translation of the computations before the
statement. In the previous section we proved that the class of the
ramification locus of $q_X$ was $2h+H+K$. By duality, the
ramification locus of $q_{X^*}$ will be $2h^*+H+K$. Multiplying
these two classes and substracting the cycle class of $R$ we
complete the proof of 1).

The degree, $\mu_1$ and class of the focal surface are easy to
obtain, by just multiplying the cycle class of $\tilde F$
respectively by $h^2$, $hh^*$ and ${h^*}^2$ (we are also using the
adjunction identity $KH+H^2=2g-2$ and the Noether formula
$c_2(T_X)=12\chi({\cal O}_X)-K^2$). To compute the other invariants we
need to know the Hilbert polynomial of $\tilde F$. For this purpose,
it is not enough to know the cycle class of $\tilde F$, but to use
the fact that it is obtained by linkage inside a complete
intersection $M$ of divisors of classes
$2h+H+K$ and $2h^*+H+K$. This fact implies (see \Peskine, Prop.
1.1) that there is an exact sequence
$$0\to {\cal I}_M\to {\cal I}_R\to {\cal H}om_{{\cal O}_{A_X}}
({\cal O}_{\tilde F},{\cal O}_M)\to 0.$$
Now the wanted invariants can be directly obtained from the coefficients
of the polynomial $\chi(\omega_{\tilde F}(Th))\in{\Bbb Q}[T]$. We will
compute it from the above exact sequence. We first observe that, by
adjunction and\canonicoincidenza,
$\omega_M\cong\omega_{A_X|M}(2h+2h^*+2H+2K)\cong{\cal O}_M(4H+3K)$, so
that
$${\cal H}om_{{\cal O}_{A_X}}({\cal O}_{\tilde F},{\cal O}_M)
\cong {\cal H}om_{{\cal O}_{A_X}}
({\cal O}_{\tilde F},\omega_M)(-4H-3K)\cong
\omega_{\tilde F}(-4H-3K).$$
We then need to compute $\chi({\cal I}_M(Th+4H+3K))$, which is very
easy since $M$ is a complete intersection. On the other hand, from
the construction of $R$, there is an exact sequence
$$0\to T_{A_X}\to q_{13}^*T_J\to {\cal I}_R(h+h^*+2H+K)\to 0$$
from where we can compute $\chi({\cal I}_R(Th+4H+3K))$. With the Maple
package Schubert one performs the computations and arrives to the
wanted result.
\qed

\bigskip

\noindent{\bf Remarks:} 1) The degree and class of the focal surface
are very well-known and there are much simpler ways to compute them.
In fact, all the other numerical invariants of the focal surface,
except $\mu_1$, can be computed by just using the incidence variety
point-line. In fact, $\mu_1$ can also be computed by using that it
is the class of $X$ considered as a surface in ${\Bbb P}^5$ (see
\Rothuno). Then $\mu_1$ is nothing but the degree of
$c_2(P^1({\cal O}_X(1))$, which is easily seen to be the value just
computed.

2) The computations previous to the proof of the above proposition
show that, for a general line $L$ of a congruence $X$, there are
exactly two pencils $\Omega(x_1,\Pi_1)$ and $\Omega(x_2,\Pi_2)$
(given by the two branch points of $q_{13}$ on $L$) that are
tangent to $X$ at the point represented by $L$. Hence the embedded
tangent plane of $X$ (as a surface in ${\Bbb P}^5$) at $L$ is the
one generated by these two pencils. However, the tangent plane at
$x_1$ of the focal surface $F$ is $\Pi_2$ and reciprocally.
\bigskip

\proclaim Proposition \singfocale. Let $X$ be a smooth congruence,
and let $F$ be its total focal surface. Assuming that the only
one-dimensional singular locus of $F$ consists of a nodal curve $D$
and a cuspidal curve $C$, then
$$\deg(D)=2a^2-10a+4b+4ag+2g^2-34g+32-4K_X^2+12\chi({\cal O}_X)$$
$$\deg(C)=3a-3b+18g-18+3K_X^2-12\chi({\cal O}_X).$$

{\it Proof:} The underlying idea is quite simple, although it
requires a precise construction of some technical complexity. We
just want to study when the fibers of the map $q_X:I_X\to{\Bbb
P}^3$ contain three infinitely close points (to find the cusps) or
two pairs of infinitely close points (to find the nodes). We will
consider more generally the projection $\pi:{\Bbb P}^3\times
X\to{\Bbb P}^3$ and apply to it a theory of infinitely close points
of its fibers (which will be just infinitely closed points in $X$).
To avoid some technical difficulties, we will reduce to the case of
cuspidal points. Note that, since we know from Prop. \classefocale\
the geometric genus of the hyperplane section of $F$, the degree of
the nodal curve can be computed at once if we know the degree of
the cupidal curve (just apply the Pl\"ucker formula to a hyperplane
section of $F$).

So we want to find a variety parametrizing sets of three infinitely
close points in the fibers of $\pi$, to find their the subset
$\tilde C$ of those who are in fact on $X$. We will follow the
construction of \Arkiv. Clearly, the variety parametrizing pairs of
infinitely close points in the fibers of $\pi$ is nothing but
${\Bbb P}(\Omega_{{\Bbb P}^3\times X/{\Bbb P}^3})={\Bbb
P}^3\times{\Bbb P}(\Omega_X)=:{\Bbb P}^3\times D^1_X$. Let
$f_1:D^1X\to X$ the structure projection and write $L_1$ for the
tautological line bundle of $D^1X$. Now the variety parametrizing
sets of three infinitely closed points in the fiber of $\pi$ is
given by ${\Bbb P}^3\times D^2X$, where $D^2X:={\Bbb P}(G)$, $G$
being the rank-two vector bundle on $D^1X$ defined as a
push-forward in the following commutative diagram of exact
sequences:
$$\leqno{\fibratoSemple}\ \ \ \ \ \ \ \ \ \ \ \matrix{
 &   &      &   &    0     &   &     0    &   & \cr
 &   &      &   &\downarrow&   &\downarrow&   & \cr
0&\to&\Omega_{D^1X/X}\otimes L_1&\to&f_1^*\Omega_X&\to&L_1&\to&0\cr
 &   &  ||  &   &\downarrow&   &\downarrow&   & \cr
0&\to&\Omega_{D^1X/X}\otimes L_1&\to&\Omega_{D^1X}&\to& G &\to&0\cr
 &   &      &   &\downarrow&   &\downarrow&   & \cr
 &   &      &   &\Omega_{D^1X/X}& = &\Omega_{D^1X/X}&   &\cr
 &   &      &   &\downarrow&   &\downarrow&   & \cr
 &   &      &   &    0     &   &     0    &   &
}$$ (see \Arkiv\ for more details). Let $f_2:D^2X\to D^1X$ denote
the structure projection and let $L_2$ be the tautological line
bundle on $D^2X$. We are now going to try to restrict the above
construction to $X$, having in mind that we are not only looking
for infinitely close points whose support is in the fiber of $q_X$:
we need the infinitesimal information defined by these points to be
also in the fiber of $q_X$.

The first step is conceptually easy. Since we want the infinitely
close points to be supported on the fiber of $q_X$, it suffices to
restrict the above construction to $I_X$ rather than working on the
whole ${\Bbb P}^3\times X$. Observe that the inclusion $I_X\subset
{\Bbb P}^3\times X$ is induced by projectivizing the quotient of
bundles in the restriction to $X$ of the universal sequence
\universale. Hence, $I_X$ is defined in ${\Bbb P}^3\times X$ as the
zero locus of the natural section of $\pi^*S_{|X}\otimes{\cal
O}_{{\Bbb P}^3}(1)$. In particular, the class of $I_X$ inside
${\Bbb P}^3\times X$ is given (we will omit to write pullbacks when
they are clear) by
$$\leqno{\classeIX}\ \ \ \ \ \ \ \ \ \ \ \ \ \ \ \ \ \ \ \ \ \ \ \ \ \
\ \ \ \ \ \ \ \ \ [I_X]=h^2+hH+c_2(\pi^*S_{|X}).$$

Keep the same notations for the above construction restricted to $I_X$
and let us see now when an element of $(1\times
f_1)^{-1}(I_X)\subset{\Bbb P}^3\times D^1X$ corresponds to a pair of
infinitely close points contained as a scheme in the fiber of
$q_X$. Those elements will be characterized by the fact that the
universal quotient $(1\times f_1)^*\Omega_{{\Bbb P}^3\times X/{\Bbb
P}^3}\to L_1$ factors through $(1\times f_1)^*\Omega_{I_X}$. This
means that the composed map
$N^*_{I_X/{\Bbb P}^3\times X}\to
(1\times f_1)^*\Omega_{{\Bbb P}^3\times X/{\Bbb P}^3}\to L_1$ is zero.

Since $I_X$ was defined as the zero locus of a section of
$\pi^*S_{|X}\otimes{\cal O}_{{\Bbb P}^3}(1)$, then its normal bundle
$N_{I_X/{\Bbb P}^3\times X}$ is isomorphic to
$q_X^*S_{|X}\otimes{\cal O}_{{\Bbb P}^3}(1)$. Hence, the
wanted subset $X'\subset (1\times f_1)^{-1}(I_X)\subset
{\Bbb P}^3\times D^1X$ is defined as the zero locus of a section of
$(1\times f_1)^*(p_X^*S_{|X}\otimes{\cal O}_{{\Bbb P}^3}(1))\otimes
L_1$, and its class inside $(1\times f_1)^{-1}(I_X)$ is then:
$$\leqno{\classeDX}\ \ \ \ \ \ \ \ \ \ \ \ \ \ \ \ \ \ \ \ \ \ \
[X']=c_1(L_1)^2+2c_1(L_1)h+c_1(L_1)H+h^2+hH+c_2(S).$$

We restrict to that subset and again abuse the notation by not changing
it after the restriction. Our final step is to identify inside
$(1\times f_2)^{-1}(X')\subset {\Bbb P}^3\times D^2X$ the subset
$X''$ of those infinitely closed points in the fiber of
$q_X$. The apparently new problem is that now $D^2X$ is not the
projectivization of a cotangent bundle, but of its quotient $G$
defined in \fibratoSemple. However this is not a problem, since the
reasoning is exactly as above. Indeed, we have now a universal
epimorphism $(1\times f_2)^*G\to L_2$ on $(1\times f_2)^{-1}(X')$, and
again $X''$ is the locus for which the natural composition
$$N^*_{X'/(1\times f_1)^{-1}(I_X)}\to\Omega_{(1\times
f_2)^{-1}(X')}\to (1\times f_2)^*G\to L_2$$
is zero. Hence $X''$ is the zero locus of a section of
$(1\times f_2)^*((1\times f_1)^*(q_X^*S_{|X}\otimes{\cal O}_{{\Bbb
P}^3}(1))\otimes L_1)\otimes L_2$, and its class in $(1\times
f_2)^{-1}(X')$ is then
$$\leqno{\classeDDX}\ \ \ \ \ \ \ \ \ \ \ \ \ \ \ \ \ \
[X'']=c_1(L_2)^2+2c_1(L_1)c_1(L_2)+2c_1(L_2)h+c_1(L_2)H$$
$$+c_1(L_1)^2+2c_1(L_1)h+c_1(L_1)H+h^2+hH+c_2(S_{|X}).$$
We finally observe that the expressions \classeIX, \classeDX\ and
\classeDDX\ can be lifted to classes in ${\Bbb P}^3\times D^2X$, so
that the degree of the cuspidal curve can be computed by intersecting
there these three classes and the class $h$ of a hyperplane in ${\Bbb
P}^3$. Now to finish the proof we use Schubert once more.
\qed

\bigskip

\noindent{\bf Remarks:} 1) If at the end of the above proof we
multiply by $H$ instead of $h$, we would get the degree of the ruled
surface consisting of those lines such that one of its two focal
points is a cusp in the focal surface. This number turns out to be
$4a+4b+12g-12$, and was already known by the classics (see
\Schumacher\ \S13 or \Rothdue\ page 197). In fact, they also
knew how to compute the degree invariants of the nodal and cuspidal
curves. Of course, they computed all these invariants in terms of
other invariants, as for instance $\mu_1$, instead of the ``modern''
invariants that we use.

2) The above proposition is valid if $X$ has not a curve of
fundamental points. This hypothesis is hidden in the statement,
since a fundamental point produces a singular point on the focal
surface whose singularity is neither a node nor a cusp. In fact,
a fundamental curve produces in the set $X''$ defined in the above
proof a component of dimension two. However, smooth congruences with
a fundamental curve are classified (see \ArrondoGross).

3) The same kind of observation can be made when we have a finite
number of fundamental points. The set $X''$ contains the cones formed
by the lines of the congruence through any fundamental point. These
cones do not count when intersecting with the class $h$, so they do
not affect to the formula of $\deg C$. However, the formula in part
1) of the remark takes also account of the sum of the degrees of
these cones.
\bigskip

Let us now apply these remarks to some examples.

\bigskip

\noindent{\bf Example \duedue}: (See also Example \Kummer\ below).
The complete intersection of $G(1,3)$ with a general hyperplane and
a general quadric produces a congruence of bidegree $(2,2)$ and
$g=1$, in particular without cuspidal curve. As a surface in ${\Bbb
P}^5$,it is the surface given by the polarized pair
$(Bl_{p1,..p5}{\Bbb P}^2, 3L-E_{1}-...-E_5),$ i.e. by the linear
system of plane cubics through five points. Therefore, the
congruence contains sixteen lines of ${\Bbb P}^5$, which correspond
to sixteen pencils of lines of ${\Bbb P}^2$. Hence the congruence
contains sixteen fundamental points (and sixteen fundamental
planes) and the degree of the corresponding cone at each of them is
one. In fact, the formula in 1) yields $16$, and hence Remark 3)
proves that there are no more fundamental points (or fundamental
planes).

\noindent{\bf Example \duetre}: As a second example, we can
consider the congruence $X$ of bidegree $(2,3)$ which is the Del
Pezzo surface (i.e. $g=1$) given by the polarized pair
$(Bl_{p1,..p4}{\Bbb P}^2, 3L-E_{1}-...-E_4),$ i.e. by the linear
system of plane cubics through four points. It is known, and easy
to verify, that such a Del Pezzo surface contains exactly ten lines
(the four exceptional lines and the six lines joining the four base
points) and five pencil of conics (the four pencils given by the
lines in ${\Bbb P}^3$ through one base point and the pencil given
by the conics through the base points). This implies that the
congruence has $15$ fundamental points and $10$ fundamental planes
in ${\Bbb P}^3.$ Indeed the ten lines give rise to ten fundamental
points and ten fundamental planes. Moreover each pencil of conics
contains at least one conic which is contained in an alpha-plane.
Let us prove this fact for instance for the pencil $|L-E_1|$ (for
the others the proof is the same). Since the image of $E_1$ in
$G(1,3)$ is a line, in particular it is contained in an
alpha-plane, so that there is a section of $S_{|X}$ vanishing on
$E_1$, i.e. a section of $S_{|X}(-E_1)$. Since
$c_2(S_{|X}(-E_1))=0$, it follows easily that there is an exact
sequence
$$0\to{\cal O}_X(E_1)\to S_{|X}\to{\cal
O}_X(3L-2E_1-E_2-E_3-E_4)\to0.$$
{}From this it follows that $h^0(S_{|X}(-L+E_1))=1$, and hence any
conic in $|L-E_1|$ is contained in the zero locus of a section of
$S_{|X}$. But observe that $h^0(S_{|X})=5$, so that exactly a
hyperplane inside $H^0(S_{|X})$ corresponds to alpha-planes. This
means that at least one section corresponding to an alpha-plane
vanishes on a a conic of the pencil, as wanted. Applying now Remark
3) we see that the degree of the ruled surface generated by the lines
throught the fundamental points is $20$ (since the bidegree is
$(2,3)$ there is no cuspidal curve). Since we have found ten cones of
degree one and five cones of degree two, there are no more
fundamental points in the congruence.

\noindent{\bf Example \tretre}
In this last example, we consider the congruence of bidegree $(3,3)$
and $g=2$ which is the rational surface given by the polarized pair
$(Bl_{p1,..p7}{\Bbb P}^2, 4L-2E_{1}-E_2-...-E_7),$ i.e. by the
linear system of plane quartics with a fixed double point and through
other six points. Such a Castelnuovo surface contains twelve lines
(the six exceptional lines corresponding to simple points, and the
six lines joining the double point with the other six ones) and $32$
conics (the one corresponding to the double base point, the $15$
corresponding to the lines joining two simple base points, the $15$
corresponding to conics through the double point and other four
simple points, and the one corresponding to the cubic with a double
point in the double base point and passing through the other base
points). Fano shows (\Fano, pages 154-155) that besides the twelve
fundamental points coming from the twelve lines of the congruence,
there can be other fundamental points (vertex of cones corresponding
to conics lying in alpha-planes) or not, depending on the projective
embedding. Specifically the Castelnuovo surface is the complete
intersection of the cubic Segre threefold and a smooth hyperquadric
in ${\Bbb P}^5$. It is hence contained in a  three-dimensional linear
system of hyperquadrics. Each smooth quadric in the system can be
viewed as a Grassmannian. While for a general quadric we do not get
extra fundamental points, for particular ones we can get one, two or
three new fundamental points.

\noindent{\bf Remark:} As shown in the previous example, the number
of the fundamental points of a congruence does not depend only on its
invariants, in particular it is meaningless to look for a formula
giving the contribution of the fundamental points only in terms of
the bidegree, of the sectional genus and of other usual invariants of
the surface. This fact seems not to be considered by Roth who gives a
formula (\Rothdue, page 198) to compute the degree $\rho_2$ of the
scroll of lines of ${\Bbb P}^3$ consisting of those lines such that
one of its two focal points is a node in the focal surface. Such a
formula, when applied to a congruence of bidegree $(a,b)$ with $a\leq
3,$ hence without nodal curve, should give the number of fundamental
points. However the formula for $\rho_2$ given by Roth fails for
several congruences (and not only for the above example).

\bigskip
\bigskip

\noindent{\semilarge \S 2. Congruence of the bisecants to a space
curve.}
\bigskip

In this section we describe the congruences of the chords of a
smooth irreducible skew curve $\Gamma$ in ${\Bbb P}^3.$

Let $\Gamma$ be a curve in ${\Bbb P}^3$ and denote by $X \subset
G(1,3)$ the congruence of the bisecants to $\Gamma.$ Throughout
this section $\Gamma$ will be assumed to be smooth, irreducible and
not contained in a plane. We will also write $d$ for the degree of
$\Gamma$ and $p$ for its genus.

It is known (see \Gross\ Theor. 2.5) that $X$ is singular unless
$\Gamma$ is a rational cubic or an elliptic quartic curve. So, from
now on, being mostly interested in the case of smooth congruences,
we could confine ourself to consider the case of these two curves,
but we prefer to study a more general situation.

\proclaim Proposition \invarianticorde. Let $\Gamma$ be as above.
Then the congruence $X$ of bisecants to $\Gamma$ has bidegree
$(a,b)=({1\over2}(d-1)(d-2)-p,{1\over2}d(d-1))$ and sectional genus
$g={1\over2}(d-2)(d-3+2p)$.

{\it Proof:} The congruence is naturally parametrized by the second
symmetric product $S=C^{(2)}$ of $C$. We will regard $S$ as the
quotient of $C\times C$ under the standard involution. Let $L$ be
the line bundle giving the embedding of $C$ into ${\Bbb P}^3$, and
write $L_1$ and $L_2$ for the corresponding pullbacks of $L$ to
$C\times C$ via the two natural projections. If $D$ denotes the
diagonal of $C\times C$, there is an epimorphism $L_1\oplus L_2\to
{\cal O}_D(L)$. Its kernel is invariant under the involution of
$C\times C$ hence it is the pullback of a rank-two vector bundle
$Q$ on $S$ (the so-called secant bundle). This vector bundle is the
one that gives the map from $S$ to $G(1,3)$ whose image is the
congruence $X$.

In the intersection ring of $S$ consider the following classes: $P$
will represent the class of pairs containing a fixed point of $C$,
and $\Delta$ will be the diagonal class, i.e. the image of $D$. We
recall the following intersection numbers: $P\cdot
P=P\cdot\Delta=1$, $\Delta\cdot\Delta=2(2-2g)$.

With this notation, the Chern classes of $Q$ are
$c_1(Q)=dP-{1\over2}\Delta$ and $c_2(Q)={1\over2}d(d-1))$. From
this one can readily obtain the bidegree by using that
$a=c_1(Q)^2-c_2(Q)$ and $b=c_2(Q)$. Notice that this bidegree could
also be obtained by simple geometric arguments.

In order to obtain the sectional genus of $X$ we need to obtain the
canonical class of $S$. This can be easily done since $C\times C$
is a double cover of $S$ ramified along the diagonal. We then have
that numerically $K_S\equiv (2-2p)P+{1\over2}\Delta$ and from here
the wanted equality for $g$ follows. \qed

\bigskip

\noindent{\bf Remark:} In the same way it is easy to find the rest
of the invariants for $S$. In particular, $K^2_S=4p^2-13p+9$ and
$\chi({\cal O}_S)={1\over2}(p-1)(p-2)$.

\bigskip

\noindent{\bf Definition:} Let $\Gamma$ be as above and consider
two distinct points $x$, $y$ of it. The chord $<x,y>$ through $x$
and $y$ is said to be {\it stationary} if the tangent lines $t_x$
and $t_y$ to $\Gamma$, at $x$ and $y$ respectively, are incident.

Denote by $T(x,y)$ the tangent plane to the congruence $X$ at the
point corresponding to a chord $<x,y>$. It is quite easy to verify
that, if the chord $<x,y>$ is stationary, then the plane $T(x,y)$
is contained in the Grassmannian $G(1,3)$, actually it is the
beta-plane generated by $t_x$ and $t_y$, as we will show in the
following (probably well known) lemmas.

\proclaim Lemma \tangenteChow. Let $\Gamma$ be as above and let $C$
be the Chow complex of lines intersecting $\Gamma$. Let $x$ be a
point of $\Gamma$, consider a line $L$ passing through $x$, denote
by $t_x$ the tangent line of $\Gamma$ at $x$, and by $\Pi$ the
plane generated by $L$ and $t_x$. Then the corresponding branch of
$C$ is smooth at the point represented by $L$ if and only if $L$ is
different from $t_x$. Moreover, in this case the embedded tangent
space of this branch of $C$ at $L$ is generated by the alpha-plane
$\alpha(x)$ and the beta-plane $\beta(\Pi)$.

{\it Proof:} This is just an easy local computation. Choose
coordinates $z_0,z_1,z_2,z_3$ in ${\Bbb P}^3$ so that the point $x$
becomes $(1:0:0:0)$ and the tangent line $t_x$ is $z_2=z_3=0$.
 Working in the open affine set $z_0=1$, we can parametrize $\Gamma$
locally at $x$ (which is now the origin) by $z_1=t, z_2=f(t),
z_3=g(t)$ with $f(0)=g(0)=f'(0)=g'(0)=0$. Let $L$ be the line passing
through $x$ and through the point of coordinates $(0:a_1:a_2:a_3)$,
and assuming $a_1\neq0$, put $a_1=1, a_2=u, a_3=v$. Then a local
parametrization in the open subset $\{p_{01}\neq0\}\subset G(1,3)$ of
the corresponding branch of $H$ at the point represented by $L$ is
given by
$$(t,u,v)\mapsto(p_{02},p_{03},p_{12},p_{13})=(u,v,ut-f(t),vt-g(t))$$
Hence, the corresponding branch of $H$ at $L$ is smooth if and only
if $(u,v)\neq(0,0)$, i.e., if and only if $L$ is different from the
tangent line at x. In this case, the embedded tangent space of $H$ at
$L$ has (affine) parametric equations

$$\left\{\eqalign{p_{01}=&1\cr
	   p_{02}=&u+\lambda\cr
	   p_{03}=&v+\mu\cr
	   p_{12}=&\nu u\cr
	   p_{13}=&\nu v\cr
	   p_{23}=&0\cr}\right.$$

\noindent i.e. it is the projective plane $vp_{12}-up_{13}=p_{23}=0$,
which is generated by the alpha-plane $\alpha(x)$ (of equations
$p_{12}=p_{13}=p_{23}=0$) and the beta-plane $\beta(\Pi)$ (of
equations $vp_{02}-up_{03}=vp_{12}-up_{13}=p_{23}=0$).

\qed

\proclaim Lemma \tangentecorde. Let $\Gamma$ be as above and $X$
the congruence of bisecants to $\Gamma$. Let $L$ be line having
exactly two intersection points $x$ and $y$ with $\Gamma$. Then $L$
represents a smooth point of $X$ if and only if it is different
from both $t_x$ and $t_y$. In this case, denote by $\Pi_x$ the
plane generated by $L$ and $t_x$ and by $\Pi_y$ the plane generated
by $L$ and $t_y$, then the embedded tangent space to $X$ at $L$ is
generated by the pencils $\Omega(x,\Pi_y)$ and $\Omega(y,\Pi_x)$.

{\it Proof:} In fact, locally at $L$ the congruence $X$ is the
complete intersection of the two branches of the Chow complex of
$\Gamma$ corresponding to the points
$x$ and $y$. Hence, $L$ is a smooth point of $X$ if and
only if the two branches are smooth at $L$ and their embedded tangent
spaces are different. This, due to Lemma \tangenteChow, happens if and only if
$L$ is neither $t_x$ nor $t_y$ and $x \neq y$.
If this is the case,
the embedded tangent plane of $X$ at $L$ will be the intersection
of the embedded tangent spaces of the two branches, which, due to Lemma
\tangenteChow, gives the thesis.

\bigskip

\noindent{\bf Remark:} The Lemma above immediately implies that, if
the chord $L=<x,y>$ is stationary, then the plane $T(x,y)$ is
contained in the Grassmannian $G(1,3)$: actually it is the
beta-plane $\beta(\Pi_x)=\beta(\Pi_y)$. Since a curve has in
general a one-dimensional family of stationary bisecants, the
corresponding congruence of bisecants will have a focal surface,
even if one would expect the focal locus to be just the curve
$\Gamma$.

From Propositions \invarianticorde\ and \classefocale, it follows
that the degree of the (total) focal surface must be
$2(d-3)(d-1+p)$, which coincides with the degree of the ruled
surface of stationary bisecants to $\Gamma$ (see \Johnsen, Remark
5.2). The twisted cubic is the only curve in ${\Bbb P}^3$ without
stationary bisecants, so we study next in detail the only other
example of smooth congruence of bisecants.

\noindent{\bf Example \duesei:} Let $X$ be the congruence of bisecants
to an elliptic quartic curve $\Gamma\subset{\Bbb P}^3$. It is then
known that $X$ is a smooth congruence of bidegree $(2,6)$ and
sectional genus $g=3$. The strict focal ``surface'' $F_0$ will be
$\Gamma$, while $F$ consists of the four quadric cones containing
$\Gamma$. Indeed it is easy to see that a bisecant to $\Gamma$ is
stationary if and only if it is contained in one of the quadric cones
containing $\Gamma$. Observe that  we then obtain the expected degree
eight for the focal surface of $X$.

\bigskip
\bigskip

\noindent {\semilarge \S 3. Congruences of bitangents and flexes to
a smooth surface in ${\Bbb P}^3$: global study}
\bigskip

\noindent Let $\Sigma\subset{\Bbb P}^3$ be a surface of degree $d$,
that we will assume, unless otherwise specified, to be smooth. In
fact we will also sometimes assume $\Sigma$ to be general enough, so
that, for $d\ge4$, its Picard group will be generated by the hyperplane
section. Following the ideas of \McCrory\ and \Welters, we consider the
projective bundle
$p:Y={\Bbb P}(\Omega_\Sigma(2))\to\Sigma$. Any point of $Y$ can be
regarded as a pair $(x,L)$, where $x$ is a point of $\Sigma$ and
$L$ is a tangent line to $\Sigma$ at $x$. Therefore there is a map
$\varphi:Y\to G(1,3)$. In fact the twist in the projective
bundle was chosen so that the tautological line bundle of $Y$ became
the pull-back of the hyperplane section of
$G(1,3)$. Let us write ${\cal O}_Y(\ell)$ for the tautological
line bundle on $Y$ and ${\cal O}_Y(h)$ for the pull-back via $p$ of
the hyperplane line bundle of $\Sigma\subset{\Bbb P}^3$. In terms of
vector bundles, the map
$\varphi$ is defined by the rank-two vector bundle $Q$ on $Y$
defined as a push-forward in the commutative diagram:

$$\leqno{\univ}\ \ \ \ \ \ \ \ \ \ \matrix{
 &   &      &   &    0     &   &     0    &   & \cr
 &   &      &   &\downarrow&   &\downarrow&   & \cr
0&\to&\Omega_{Y/\Sigma}(\ell-h)&\to&p^*\Omega_\Sigma(h)&\to&{\cal
O}_Y(\ell-h)&\to&0\cr
 &   &  ||  &   &\downarrow&   &\downarrow&   & \cr
0&\to&\Omega_{Y/\Sigma}(\ell-h)&\to&p^*(P^1({\cal O}_\Sigma(1)))
&\to&     Q    &\to&0\cr
 &   &      &   &\downarrow&   &\downarrow&   & \cr
 &   &      &   &{\cal O}_Y(h)& = &{\cal O}_Y(h)&   & \cr
 &   &      &   &\downarrow&   &\downarrow&   & \cr
 &   &      &   &    0     &   &     0    &   &
}$$

\noindent Here the top horizontal sequence is the universal sequence
on the projective bundle $Y$ tensored with ${\cal O}_Y(-h)$ and the
middle vertical sequence is the pull-back of the one defining the
bundle of principal parts of ${\cal O}_\Sigma(1)$. The map $\varphi$
is precisely defined by the composed epimorphism $H^0({\Bbb P}^3,
{\cal O}_{{\Bbb P}^3}(1))\otimes{\cal O}_Y\to p^*(P^1({\cal
O}_\Sigma(1)))\to Q$. The following closed surfaces of $Y$ will play an
important role in the sequel:

$$Y':=\{(x,L)\in Y\ \mid\ x \ {\rm is\ a\ parabolic\ point\ of\
\Sigma}\}$$
$$Y_1:=\{(x,L)\in Y\ \mid\ L \ {\rm is\ a\ bitangent\ line\ of\
\Sigma}\}$$
$$Y_2:=\{(x,L)\in Y\ \mid\ L \ {\rm is\ an\ inflection\ line\ of\
\Sigma}\}.$$

\noindent Of course, all the above sets are defined as a closure
(for the definition of parabolic point, see for instance \McCrory).

\bigskip

\proclaim Proposition \classi. The classes of $Y', Y_1, Y_2$ in
the Picard group of $Y$ are: $[Y']=4(d-2)h$,
$[Y_1]=(d+2)(d-3)\ell-4(d-3)h$ and
$[Y_2]=2\ell+(d-4)h$.

{\it Proof:} The surface $Y'$ is just the pullback via $p$ of the
parabolic curve on $\Sigma$, where it has class $4(d-2)h$, as shown
in \McCrory\ (anyway, the idea is that the parabolic curve is defined
by the Hessian matrix to be singular).

The class of $Y_1$ is computed in \Welters\ Prop. 3.14 for $d=4$. We
essentially reproduce here Welters' ideas. Since its class is
not so crutial, we chose the simplest but least general of his proofs.

If $\Sigma$ is sufficiently general and $d\ge 4$, then the Picard group
of $Y$ is generated by the classes of $\ell$ and
$h$. Hence the class of $Y_1$ will be of the form $m\ell+nh$. The
first integer $m$ is in fact the degree of the projection
$Y_1\to\Sigma$, hence it is the number of tangents at a general
point of $x\in\Sigma$ that are tangent to $\Sigma$ at another point.

To compute this number, consider $\Pi$ the tangent plane to
$\Sigma$ at $x$ and let $C$ be the intersection of $\Sigma$ with
$\Pi$. Hence $C$ is a plane curve of degree $d$ with one ordinary
node at $x$ (hence of geometric genus ${d(d-3)\over2}$) and
$m$ is the number of lines which are tangent to $C$ outside
$x$ and pass through $x$. In other words, $m$ is the number
of branch points of the $(d-2):1$ morphism $C\to{\Bbb P}^1$
defined by the projection from $x$. From Hurwitz theorem one
immediately gets $m=(d+2)(d-3)$.

To compute $n$ we can use the fact that the order of the congruence
of bitangents to $\Sigma$ is $1/2d(d-2)(d-3)(d+3)$ (the number of
bitangents of a general plane curve of degree $d$). Since the map
$Y_1\to G(1,3)$ (restriction of $\varphi$) is a double cover of
such a congruence, it follows that $c_2(Q_{Y_1})=d(d-2)(d-3)(d+3)$.
This Chern class can be computed (with the help of the Maple
package Schubert) from diagram \univ\ in terms of $n$, and making it
equal to the second term one gets the required value of $n$.

The class of $Y_2$ can be computed in a more direct way. First we
recall that inflectional tangent vectors to $\Sigma$ are those in
the kernel of the second fundamental form $II:Symm^2
T_\Sigma\to N_\Sigma$. Here, $N_\Sigma={\cal O}_\Sigma(d)$ is the
normal bundle of $\Sigma$. Hence we are looking at the points of
$\Sigma$ for which the section of $Symm^2(\Omega_\Sigma)(d)$
corresponding to $II$ is zero. Regarding this section as a section
of $Symm^2(\Omega_\Sigma(2))(d-4)$, we see from the projection
formula that this corresponds to a section of ${\cal
O}_Y(2\ell+(d-4)h)$, whose zero locus is precisely $Y_2$.
\qed

\bigskip

Let us write $X_i=\varphi(Y_i)$ and $\varphi_i=\varphi_{|Y_i}$ for
$i=1,2$. Then $X_1$ is the congruence of bitangents of $\Sigma$ and
$X_2$ is the congruence of inflectional lines of $\Sigma$. They
both are contained in the complex ${\cal H}:=\varphi(Y)$ of lines
tangent to $\Sigma$. What makes this approach so different among
these two congruences is that, while the map $\varphi_1:Y_1\to X_1$
is a double cover, the map $\varphi_2:Y_2\to X_2$ is birational (in
both cases, the map $\varphi_i$ is finite as long as $\Sigma$ does
not contain any line). Hence we can easily compute the bidegree of
both congruences, but it will be possible only for $X_2$ to compute
all its invariants. As remarked in \Welters\ page 30, the map
$\varphi_1$ is branched over the curve of hyperflexes; the study of
such a curve would certainly allow to compute all the invariants of
$X_1$ from the ones of $Y_1$. We will use however a different way
(see Proposition \bitangenti\ below), which is more elegant and
will also allow us to remove the genericity hypothesis for
$\Sigma$.

\proclaim Proposition \invarianti. The congruence $X_1$ has bidegree
$({1\over 2}d(d-1)(d-2)(d-3),{1\over 2}d(d-2)(d-3)(d+3)$, while the bidegree
 of $X_2$ is $(d(d-1)(d-2),3d(d-2))$ and the
sectional (geometric) genus $g=5d^3-18d^2+14d+1$.
 Moreover, the congruence $X_2$ is
never smooth.

{\it Proof:} The map $\varphi_i$, as a map to $G(1,3)$, is given by
the rank-two vector bundle $Q_{Y_i}$. Since $\varphi_1$ is a double
cover, then the class of $X_1$ is ${1\over
2}c_2(Q_{Y_1})$, and in fact we have already seen (or rather
impose) in the proof of Prop. \classi\ that this is
${1\over 2}d(d-2)(d-3)(d+3)$. Analogously, its order is ${1\over
2}(c_1(Q_{Y_1})^2-c_2(Q_{Y_1}))={1\over 2}d(d-1)(d-2)(d-3)$, as
easily computed again with the help of the Maple package Schubert.

In a similar but easier way, since $\varphi_2$ is now birational, the
class of $X_2$ is just $b=c_2(Q_{Y_2})=3d(d-2)$ (which in fact
corresponds to the number of flexes of a general plane curve of
degree $d$) while its order is
$a=c_1(Q_{Y_2})^2-c_2(Q_{Y_2})=d(d-1)(d-2)$. On the other hand,
assume now that $X_2$ is smooth. Hence, if $\Sigma$ does not
contain any line, the map $\varphi_2:Y_2\to G(1,3)$ is necessarily
an immersion, and the double-point formula for it would yield
$a^2+b^2-c_2(N)=0$, where $N$ is the cokernel of the bundle inclusion
$T_{Y_2}\to \varphi^*T_{G(1,3)}$. But taking into account that
$\varphi^*T_{G(1,3)}\cong S_{Y_2}\otimes Q_{Y_2}$ (where $S_{Y_2}$ is
the dual of the kernel of a natural epimorphism ${\cal
O}_{Y_2}^{\oplus 4}\to Q_{Y_2}$, and in fact the pull-back to $Y_2$
of the universal bundle $S$ on $G(1,3)$), with the help once more of
the Schubert package we get that the double-point formula reads
$$d(d-3)(d^4-3d^3+13d^2-48d+40)=0$$
\noindent which is absurd if $d\neq 3$. The case $d=3$ (or more
generally when $\Sigma$ contains a line) is treated separately in the
following lemma. All the invariants of $X_2$ (in particular the
sectional genus) are computed using the isomorphism with $Y_2$ and
the fact that $Y_2$ is a smooth divisor of $Y$ of a known class.
\qed

\bigskip

\proclaim Lemma \retta. If $L$ is a line contained in $\Sigma$,
then the corresponding point of $X_2$ is singular of multiplicity
$3(d-2)$.
\bigskip

{\it Proof:} Let us consider the point $p_L\in X_2$ corresponding to
the line $L\subset\Sigma$. By abuse of notation, let us still call
$L$ to $\varphi_2^{-1}(p_L)$. In other words, we are identifying
$L$ with the curve in $Y_2$ which contracts to $p_L$. Since
$\varphi_2$ is birational, the multiplicity of $p_L$ will be
precisely minus the self-intersection of $L$ in $Y_2$. By adjunction
we have $L^2+K_{Y_2}L=-2$, so it is enough to prove that
$K_{Y_2}L=3d-8$. But this is an immediate consequence of the
equality $K_{Y_2}=c_1(\Omega_{Y|Y_2})-(2\ell+(d-4)h)=(3d-8)h_{|Y_2}$,
since we can then compute $K_{Y_2}L$ as the intersection in $Y$ of
$(3d-8)h$ with $L$.
\qed

\bigskip

\noindent{\bf Remark:} A similar statement was proved in \Welters\
(1.1) and (1.2) for the congruence $X_1$ of bitangents in case
$d=4$.

\bigskip

\proclaim Proposition \bitangenti. The congruence $X_1$ of
bitangents to a smooth surface $\Sigma\subset{\Bbb P}^3$ of degree
$d$ is smooth only for $d=4$. The geometric genus of its hyperplane
section is $g=d^5-{5\over2}d^4-{35\over2}d^3+60d^2-36d+1$.

{\it Proof:} The idea is to work on the Hilbert scheme
$T=Hilb^2{\Bbb P}^3$ parametrizing (unordered) couples of points of
${\Bbb P}^3$ (and then study the subset of those that produce a
bitangent line to $\Sigma$). Since two points (possibly infinitely
close) determine a line, there is a map $q:T\to G(1,3)$. On the
other hand, the set of pairs of points on a fixed line is a ${\Bbb
P}^2$, parametrized by the quadratic forms (up to a constant) on
the line. Therefore, the map $q$ endows $T$ with a projective
bundle structure $T={\Bbb P}(Symm^2Q^*)$. In this projective bundle
we have the universal quadratic form given by the bundle inclusion
$${\cal O}_T(-1)\hookrightarrow q^*Symm^2Q$$
\noindent which assigns at each couple of points the quadratic form
(defined on the line spanned by them) vanishing on those points. We
can similarly construct from this a bundle inclusion
$${\cal O}_T(-2)\hookrightarrow q^*Symm^4Q$$
\noindent which corresponds for every couple to the quartic forms
vanishing doubly at each of the points of the couple. The
multiplication of $(d-4)$-forms by this universal form determines
then another bundle inclusion $i$ which defines the bundle $R$ as a
cokernel:
$$0\to q^*Symm^{d-4}Q\otimes{\cal O}_T(-2)\mapright{i}q^*Symm^dQ\to
R\to0$$

A surface $\Sigma\subset{\Bbb P}^3$ of degree $d$ corresponds to a
section ${\cal O}_{G(1,3)}\to Symm^dQ$, and we are interested in
the locus at which the pull-back of this section lies in the image
of $i$. In other words, the zero locus of the corresponding section
of $R$ (obtained as the composition ${\cal O}_T\to q^*Symm^dQ\to R$)
is the set $\tilde X_1$ of couples of points of $\Sigma$ such that
the line defined by them is tangent at those points. The congruence
$X_1$ is the image by $q$ of $\tilde X_1$. If $X_1$ is smooth (and
$\Sigma$ does not contain any line), then $p$ defines in fact an
isomorphism between $\tilde X_1$ and $X_1$, so everything reduces to
computing the invariants of $\tilde X_1$. This is easily done by
using that $\tilde X_1$ is defined as the zero locus of the rank-four
vector bundle $R$, of which we can compute its Chern classes from the
exact sequence defining it.

To be honest, there is a technical problem that cannot be completely
solved by using the package Schubert: the Chern classes of a
symmetric power of a bundle can be computed only for a fixed
exponent, but not depending on a parameter $d$. We write the exact
result we need in Lemma \classiChern\ below, so that the interested
reader can reproduce from it all our calculations. These
calculations will provide easily the sectional genus (from the
product of the canonical class of
$\tilde X_1$ and the pull-back of the hyperplane section of
$G(1,3)$), as well as the rest of the invariants. In particular, one gets that,
if $N$ is the normal bundle of $X_1$ in $G(1,3)$, then
$a^2+b^2-c_2(N)={1\over2}d(d-4)(d^6-4d^5+2d^4-20d^3+9d^2+396d-540)$.
Hence $X_1$ is only smooth for $d=4$.
\qed

\bigskip

\proclaim Lemma \classiChern. Let $Q$ be a rank-two vector bundle
on a smooth variety and let $c_1,c_2$ be its Chern classes. Then a
symmetric power of $Q$ has Chern classes:
$$\eqalign{c_1(Symm^dQ)=&{1\over2}d(d+1)c_1\cr
	   c_2(Symm^dQ)=&{1\over24}d(d-1)(d+1)(3d+2)c_1^2+
 {1\over6}d(d+1)(d+2)c_2\cr
	   c_3(Symm^dQ)=&{1\over48}d^2(d-1)(d-2)(d+1)^2c_1^3+
 {1\over12}d^2(d-1)(d+2)(d+1)c_1c_2\cr
	   c_4(Symm^dQ)=&
{1\over1570}d(d-1)(d-2)(d-3)(d+1)(15d^3+15d^2-10d-8)c_1^4\cr
&+{1\over720}d(d-1)(d-2)(d+2)(d+1)(15d^2-5d-12)c_1^2c_2\cr
&+{1\over 360}d(d-1)(d-2)(d+1)(5d+12)c_2^2.}$$

{\it Proof:} This is just a straightforward (but terribly annoying)
calculation using the splitting principle.
\qed

\bigskip
\bigskip

\noindent {\semilarge \S 4. Congruences of bitangents and flexes to
a smooth surface in ${\Bbb P}^3$: local study}
\bigskip

In this section we analyze when a bitangent or inflectional line to
a surface becomes a focal line of the corresponding congruence. We
will find then that in both types of congruences we always get at
least one component of the focal surface made out of focal lines.
On the other hand, we will observe that the surface $\Sigma$ will
have a big multiplicity as a component of the focal surface. We
will finally check that these two atypical situations are reflected
in the formula for the degree of the focal surface, which can be
derived from the invariants of the congruences computed in the
previous section.

We prove first a series of local results about tangent spaces that
will be useful later on.

\bigskip

\proclaim Lemma \tangentecomplesso. Let $\Sigma$ be a surface in
${\Bbb P}^3$ and let ${\cal H}$ be the complex of lines tangent to
$\Sigma$.
\item{a)} If $x$ is a smooth point of $\Sigma$, $\Pi=T_x\Sigma$ the
tangent plane of $\Sigma$ at $x$, and $L$ a line contained in $\Pi$
passing through $x$, then the corresponding branch of ${\cal H}$ is
smooth at the point represented by $L$ if and only if the
intersection multiplicity at $x$ of $L$ and $\Sigma$ is exactly
two. Moreover, in this case the embedded tangent space of this
branch of ${\cal H}$ at $L$ is generated by the alpha-plane
$\alpha(x)$ and the beta-plane $\beta(\Pi)$.
\item{b)} The surface $Y_2$ is singular at the points $(x,L)$ for
which $x$ is a parabolic point (and hence $L$ is the unique
asymptotic line at $x$).

{\it Proof:} This is just based on a tedious local computation to
study the differential of $\varphi$ at the point $(x,L)$. Choose
coordinates $z_0,z_1,z_2,z_3$ in ${\Bbb P}^3$ so that the point $x$
becomes $(1:0:0:0)$, the plane $\Pi$ has equation $z_3=0$ and the
line $L$ is $z_2=z_3=0$. Working in the open affine set
$\{z_0=1\}$, we can parametrize $\Sigma$ locally at $x$ (which is
now the origin) by $z_3=f(z_1,z_2)$. Hence a local parametrization
of the corresponding branch of ${\cal H}$ at the point represented
by $L$ is given by assigning to local parameters $\lambda,u,v$ the
line generated by the rows of the matrix
$$\leqno{\matrice}\ \ \ \ \ \ \ \ \ \ \ \ \ \ \ \ \ \ \ \ \ \
\ \ \ \ \ \ \ \ \ \ \ \ \ \ \ \ \pmatrix{
1&u&    v   &        f       \cr
0&1&\lambda & f_u+\lambda f_v}$$
\noindent ($f_u$ and $f_v$ denoting the partial derivatives of f
with respect to $u$ and $v$ respectively). In this way, the
Pl\"ucker coordinates of this line in the open affine set of
$G(1,3)$ given by $\{p_{01}=1\}$ are:
$$\eqalign{p_{02}=&\lambda\cr
	   p_{03}=&f_u+\lambda f_v\cr
	   p_{12}=&\lambda u-v\cr
	   p_{13}=&uf_u+\lambda uf_v-f}$$
\noindent (These are therefore local equations for $\varphi$ at
$(x,L)$).  The Jacobian matrix with respect to $\lambda,u,v$ is then
$$\pmatrix{1&         f_u         &   u   & uf_v\cr
	   0&f_{uu}+\lambda f_{uv}&\lambda&f_u+uf_{uu}+\lambda
f_v+\lambda f_{uv}-f_u\cr
	   0&f_{uv}+\lambda f_{vv}&   -1
&uf_{uv}+\lambda uf_{vv}-f_v}$$
\noindent We now specialize to the point represented by $L$ (i.e.
$\lambda=u=v=0$) taking into account that $f_u(0,0)=f_v(0,0)=0$
(since $z_3=0$ is the tangent plane at
$p$) and get the matrix
$$\pmatrix{1&    0      & 0 &0\cr
	   0&f_{uu}(0,0)& 0 &0\cr
	   0&f_{uv}(0,0)&-1 &0}$$
\noindent Hence, the corresponding branch of ${\cal H}$ at $L$ is
smooth if and only if $f_{uu}(0,0)\neq0$, which is clearly
equivalent to the fact that $L$ meets $\Sigma$ with multiplicity
exactly two. In this case, the tangent space of ${\cal H}$ at $L$
(in the embedded tangent space of $G(1,3)$ at $L$, which is
$p_{23}=0$) has equation $p_{13}=0$. Hence the embedded tangent
space of ${\cal H}$ at $L$ is $p_{13}=p_{23}=0$, which is generated
by the alpha-plane $\alpha(x)$ (of equations
$p_{12}=p_{13}=p_{23}=0$) and the beta-plane $\beta(\Pi)$ (of
equations $p_{03}=p_{13}=p_{23}=0$). This proves a)

As for b), with the same coordinates as above, the equation of $Y_2$
is $f_{uu}^2+2f_{uv}\lambda+f_{vv}\lambda^2=0$. If $x$ is parabolic
and $L$ is the unique asymptotic line at $x$, then
$f_{uu}(0,0)=f_{uv}(0,0)=0$. Hence, the equation of $Y_2$ does not
have linear monomials and therefore the point $(x,L)$ is singular.
\qed

\bigskip

\proclaim Lemma \tangentebitangenti. Let $\Sigma$ be a surface in
${\Bbb P}^3$ and $X_1$ the congruence of bitangents to $\Sigma$.
\item{a)} Let $L$ be line having exactly two tangency points $x$ and
$y$ with $\Sigma$ ($x$ and $y$ being smooth). Then $L$ represents a
smooth  point of $X_1$ if and only if the intersection multiplicity of
$L$ and $\Sigma$ at both $x$ and $y$ is two. In this case, the
embedded tangent space to $X_1$ at $L$ is generated by the pencils
$\Omega(x,T_y\Sigma)$ and $\Omega(y,T_x\Sigma)$.
\item{b)} Let $L$ be a line of $X_1$ having only one tangency point
with $\Sigma$. Then $L$ an $\Sigma$ has intersection multiplicity at
least four at the contact point. Moreover, if the  intersection
multiplicity is exactly four, then the line $L$ is a smooth point of
$X_1$ and is not contained in the focal surface.

{\it Proof:} To prove a), we first observe that $L$ represents to a
double point of the complex of tangents ${\cal H}$, whose branches
correspond to the image by $\varphi$ of the points $(x,L)$ and
$(y,L)$. In fact, locally at $L$ the congruence $X_1$ is the
complete intersection of these two branches. Hence, $L$ will be a
smooth point of $X_1$ if and only if the two branches are smooth at
$L$ and their embedded tangent spaces are different. This second
statement is always true since $x\neq y$. Therefore, by Lemma
\tangentecomplesso, $L$ is smooth if and only if the intersection
multiplicity of $L$ and $\Sigma$ at both $x$ and $y$ is two. If
this is the case, the embedded tangent plane of $X_1$ at $L$ will
be the intersection of the embedded tangent spaces of the two
branches. Using again Lemma \tangentecomplesso\ and the fact that
$x\neq y$, the intersection with $G(1,3)$ of the embedded tangent
spaces of the two branches of ${\cal H}$, which are
$\alpha(x)\cup\beta(T_x\Sigma)$ and $\alpha(y)\cup\beta(T_y\Sigma)$
is either $\Omega(x,T_y\Sigma)\cup \Omega(y,T_x\Sigma)$ (if
$T_x\Sigma\neq T_y\Sigma$) or $\beta(T_x\Sigma)$ (if
$T_x\Sigma=T_y\Sigma$). In either case, a) follows. This fact could
also be deduced from the second remark after Prop. \classefocale.

As for b), let $L$ be a line of the congruence with only one tangency
point $x$ with $\sigma$. From the bundle construction in the proof of
Proposition \bitangenti, the equation of $\Sigma$ restricted to $L$
is divisible by four times the equation of $x$ (the universal
quadratic form on $L$ is the form vanishing twice at $x$, so that its
square vanishes four times). Hence the intersection multiplicity of
$L$ and $\Sigma$ is at least four at $x$.

Now we assume that $L$ and $\Sigma$ have intersection multiplicity
four at $x$ and choose coordinates as in the proof of Lemma
\tangentecomplesso\ (so that $x$ has affine coordinates $(0,0,0)$,
$L$ is the line $z_2=z_3=0$ and the tangent plane of $\Sigma$ at
$x$ is $z_3=0$). We can assume the local equation of $\Sigma$ at
$x$ is
$$z_3=f(z_1,z_2)=a_1z_1z_2+a_2z_2^2+a_3z_1^2z_2+a_4z_1z_2^2+a_5z_2^3$$
$$+z_1^4+a_6z_1^3z_2+a_7z_1^2z_2^2+a_8z_1z_2^3+a_9z_2^4+\ldots$$
\noindent The line of affine Pl\"ucker coordinates
$p_{02},p_{03},p_{12},p_{13}$ is then the one of affine equations
$$\eqalign{
z_2=&-p_{12}+p_{02}z_1\cr
z_3=&-p_{13}+p_{03}z_1}$$
That line will be in the congruence $X_1$ if and only if the above
substitution in the polynomial $P(z_1,z_2,z_3)=-z_3+f(z_1,z_2)$ has
two double roots. But we now observe that
$$P(z_1,p_{02}+p_{03}z_1,p_{12}+p_{13}z_1)=$$
$$=(p_{13}+a_2p_{12}^2-a_5p_{12}^3+a_9p_{12}^4)$$
$$+(-p_{03}-a_1p_{12}+a_4p_{12}^2-2a_2p_{02}p_{12}+
3a_5p_{02}p_{12}^2-a_8p_{12}^2-4a_9p_{02}p_{12}^3)z_1$$
$$+(a_1p_{02}-a_3p_{12}+a_2p_{02}^2+a_7p_{12}^2-2a_4p_{02}p_{12}
+6a_9p_{02}^2p_{12}^2-3a_5p_{02}^2p_{12}+3a_8p_{02}p_{12}^2)z_1^2$$
$$+(a_3p_{02}-a_6p_{12}+a_4p_{02}^2-2a_7p_{02}p_{12}+a_5p_{02}^3
-3a_8p_{02}^2p_{12}-4a_9p_{02}^3p_{12})z_1^3$$
$$+(1+a_6p_{02}+a_7p_{02}^2+a_8p_{02}^3+a_9p_{02}^4)z_1^4+\ldots$$

The main point now is the technical Lemma \tangentepolinomi, which we
state and prove after the end of this proof. That technical lemma
implies that $X_1$ is defined locally at $L$ by two polynomials whose
linear parts are $p_{13}$ and $-p_{03}-a_1p_{12}$. Hence, $X_1$ is
smooth at $L$, and the embedded tangent space at that point is
$p_{13}=p_{23}=p_{03}+a_1p_{12}=0$, which clearly is not contained in
$G(1,3)$. (This tangent plane can be viewed as the only plane in the
pencil determined by $\alpha(x)$ and $\beta(T_x\Sigma)$ which contains
the infinitely close line to $L$ in the quadric
$z_3=a_1z_1z_2+a_2z_2^2$, which is the osculating quadric to $\Sigma$
at $x$).
\qed

\proclaim Lemma \tangentepolinomi. Let $A_d$ be the projective
space of nonzero polynomials (up to multiplication by a nonzero
constant) in $K[T]$ of degree at most $d$ (for a fixed $d\ge4$) and
let $B_d$ be the subset of polynomials with a factor of degree four
which is a perfect square. Let the coordinates $(b_0:\ldots:b_d)$
define the polynomial $b_0+\ldots+b_dX^d\in A_d$. Then, locally at
a polynomial $X^4P$ (with $P(0)\neq0$ and $P$ square-free), $B_d$
is defined by two affine equations in $K[b_0,\ldots,\hat
b_4\ldots,b_d]$ whose linear parts are $b_0$ and $b_1$.

{\it Proof:} We start with the easy case in which $d=4$. Then we can
work in the affine space of monic polynomials, and the polynomial
$b_0+b_1X+b_2X^2+b_3X^3+X^4$ is in $B$ if and only if it the
square of a polynomial $c_0+c_1X+X^2$. Therefore one gets the
relations:
$$\eqalign{
b_0=&c_0^2\cr
b_1=&2c_0c_1\cr
b_2=&2c_0+c_1^2\cr
b_3=&2c_1}$$
\noindent From the last two equations one can obtain $c_0$ and $c_1$
as polynomials in $b_0,b_1$ without constant term, and substituting in
the first two equations one gets the wanted local equations, with
linear terms $b_0$ and $b_1$.

For a general $d$, we consider the obvious multiplication map
$$\psi:A_4\times A_{d-4}\to A_d.$$
\noindent A polynomial as in the statement is the image of a
$(X^4,P)$, and we can assume $P$ to have constant term equal to
$1$. As before, we take the obvious affine coordinates in $c_1,\ldots
c_{d-4}$ in $A_{d-4}$ and $d_0,d_1,d_2,d_3$ near $P$ and $X^4$.
Observe that $X^4$ becomes the origin, but $P$ can have arbitrary
coordinates $c_{10},\ldots,c_{d-4,0}$. The map $\psi$ is defined
in these affine sets by:
$$\eqalign{
b_0=&d_0\cr
b_1=    &c_1d_0+d_1\cr
b_2=    &c_2d_0+c_1d_1+d_2\cr
b_3=    &c_3d_0+c_2d_1+c_1d_2+d_3\cr
b_4=    &c_4d_0+c_3d_1+c_2d_2+c_1d_3+1\cr
b_5=    &c_5d_0+c_4d_1+c_3d_2+c_2d_2+c_1\cr
	&\ldots\cr
b_{d-4}=&c_{d-4}d_0+c_{d-5}d_1+c_{d-6}d_2+c_{d-7}d_3+c_{d-8}\cr
b_{d-3}=&c_{d-4}d_1+c_{d-5}d_2+c_{d-6}d_1+c_{d-7}\cr
b_{d-2}=&c_{d-4}d_2+c_{d-5}d_3+c_{d-6}\cr
b_{d-1}=&c_{d-4}d_3+c_{d-5}\cr b_d= &c_{d-4}}$$ \noindent We can
work on the open affine $b_4=1$ and divide the rest of the
coordinates by the above expression for $b_4$. It is not difficult
to check that the Jacobian matrix of $\psi$ with respect to
$d_0,d_1,d_2,d_3,c_1,\ldots,c_{d-4}$ at
$(0,0,0,0,c_{10},\ldots,c_{d-4,0})$ is lower triangular with $1$'s
in the diagonal (just observe that dividing by $b_4$ does not
change that much the aspect of the matrix). Therefore $\psi$ is
locally an isomorphism. Our hypothesis implies that $(X^4,P)$ is
the only element of $B_4\times A_{d-4}$ whose image is $X^4P$. We
therefore get a local isomorphism between $B_4\times A_{d-4}$ and
$B_d$. From what we already proved for $d=4$, the tangent space of
$B_4$ at $X^4$ is given by $d_0=d_1=0$. Looking at the differential
of $\psi$ we then conclude that the tangent space of $B_d$ at
$X^4P$ is defined by $b_0=b_1=0$, as wanted. \qed

\bigskip

\proclaim Corollary \singbitangenti. If $d\ge5$, the congruence $X_1$
has a singular curve consisting of bitangent lines to $\Sigma$ having
multiplicity three at one of the tangency points. The degree of this
curve in $G(1,3)$ is $d(d-3)(d-4)(d^2+6d-4)$.

{\it Proof:} The first statement follows at once from Lemma
\tangentebitangenti. The degree of the curve can be found, for
instance, in \Salmon, art. 598 (pages 286-287). To see a modern proof,
 a simple way
would be the following. Observe that a pair $(x,L)\in Y_1$ will
belong also to $Y_2$ if and only if either the multiplicity of
intersection of $L$ and $\Sigma$ at $x$ is at least three (when
there is another tangency point) or the intersection multiplicity
is at least four (when there is only one tangency point). The
second possibility produces a curve, whose degree in $G(1,3)$ is
given in Corollary \singflessi\ below. Once this degree is
subtracted from the intersection $[Y_1][Y_2]\ell$ in $Y$, the
remaining degree is $2d(d-3)(d-4)(d^2+6d-4)$. But, as the proof of
Lemma \tangentecomplesso\ shows, the point of $Y_2$ are in the
ramification locus of $\varphi:Y\to G(1,3)$, so that the above
degree is counted twice. \qed

\bigskip

\proclaim Lemma \tangenteflessi. Let $\Sigma$ be a surface in ${\Bbb
P}^3$ and $X_2$ the congruence of inflectional lines to $\Sigma$.
\item{a)} If $L$ is an inflectional line to $\Sigma$ at a
non-parabolic point $x$, then $L$ represents a smooth point of $X_2$
if and only if the intersection multiplicity of $L$ and $\Sigma$ at
$x$ is exactly three. In this situation, $L$ is never contained in
the focal locus of $X_2$, and the ramification index of
$I_{X_2}\to{\Bbb P}^3$ at $(x,L)$ is two.
\item{b)} If $L$ is an inflectional line to $\Sigma$ at a parabolic
point $x$ and the intersection multiplicity of $L$ and $\Sigma$ at
$x$ is exactly three, then $L$ represents a smooth point of $X_2$ and
the embedded tangent plane of $X_2$ at $L$ is the beta-plane
$\beta(T_x\Sigma)$.

{\it Proof:}
In order to prove a), let as choose coordinates as in Lemma
\tangentecomplesso. Since $x$ is not parabolic, we can also assume that
the other asymptotic line of
$\Sigma$ at $x=(1:0:0:0)$ is $z_1+z_2=z_3=0$ (this apparently
strange choice is made in order to guarantee that ${1\over f_{vv}}$
below has a Taylor expansion). In other words, there is a local affine
parametrization of $\Sigma$ at $x$ given by
$$(z_1,z_2,z_3)=(u,v,f(u,v))=
(u,v,uv+v^2+a_0u^3+a_1u^2v+a_2uv^2+a_3v^3+\ldots)$$
\noindent (where $+\ldots$ means that we are omitting terms of higher
degree). The asymptotic lines at a point parametrized by $(u,v)$
are given as the span of the rows of matrix \matrice, where
$\lambda$ is one of the roots of the equation
$f_{uu}+2f_{uv}\lambda+f_{vv}\lambda^2=0$. Taking into account that
$$\eqalign{
f_{uu}=&6a_0u+2a_1v+\ldots\cr
f_{uv}=&1+2a_1u+2a_2v+\ldots\cr
f_{vv}=&2+2a_2u+6a_3v+\ldots}$$
\noindent and using the Taylor expressions
$\sqrt{1+z}=1+{1\over2}z+\ldots$ and
${1\over2+z}={1\over2}-{1\over4}z+\ldots...$ to find a determination
of $\lambda$ in the above equation, one finds that the asymptotic
lines are locally parametrized by the rows of the matrix:
$$\pmatrix{
1&u&         v        &uv+v^2+a_0u^3+a_1u^2v+a_2uv^2+a_3v^3+\ldots\cr
0&1&-3a_0u-a_1v+\ldots&v+(a_1-6a_0)uv+(a_2-2a_1)v^2+\ldots}.$$
\noindent This gives a local affine parametrization of $X_2$:
$$\eqalign{
p_{02}=&-3a_0u+a_1v+\ldots\cr
p_{03}=&v+(a_1-6a_0)uv+(a_2-2a_1)v^2+\ldots\cr
p_{12}=&-v-3a_0u^2-a_1uv+\ldots\cr
p_{13}=&-v^2-6a_0u^3-6a_0u^2v-2a_1uv^2-a_3v^3+\ldots}$$
which must be an isomorphism at smooth points of $X_2$. Hence,
looking at the linear part, $L$ represents a smooth point if and only
if $a_0\neq0$, i.e. if and only if the line $L$ does not meet
$\Sigma$ with multiplicity greater than or equal to four. In this
case, the embedded tangent plane is then
$p_{03}+p_{12}=p_{13}=p_{23}=0$, which is not contained in $G(1,3)$.
(This tangent plane can be interpreted as at the end of the proof of
Lemma \tangentebitangenti).

To compute the ramification index of $I_{X_2}\to{\Bbb P}^3$, at
$(x,L)$, just observe that the alpha plane $\alpha(x)$ is given, in
the above local coordinates of $G(1,3)$, by the equations
$p_{12}=p_{13}=0$. Look at the above value of these coordinates in the
local parametrization of $X_2$ and using that $a_0\neq0$, we obtain a
curvilinear scheme of degree three supported at $(x,L)$. Therefore,
the ramification index is two. This completes the proof of a).

Statement b) is proved in a similar way, but observing now that,
since we are in the ramification locus of $p_{|Y_2}$, $u,v$ is not a
system of parameters for $X_2$ at $L$. Anyway, take coordinates as in
Lemma \tangentecomplesso\ or a), and we can assume that our $f$ takes
now the form
$$f(u,v)=v^2+a_0u^3+a_1u^2v+a_2uv^2+a_3v^3+\ldots$$
\noindent The new coordinate we have to choose now will be $w$, where
$$w^2=f_{uv}^2-f_{uu}f{vv}=-12a_0u-4a_1v+\ldots$$
\noindent Since by hypothesis $a_0\neq0$, we can take $v,w$ as a
system of parameters and substitute $u=-{a_1\over3a_0}v+\ldots$ in
$f,f_u,f_v,f_{uu},f_{uv},f{vv}$. In particular, we get
$$\lambda={-f_{uv}+w\over f_{vv}}={a_1^2-3a_0a_2\over
a_0}v+{1\over2}w+\ldots$$
We get now a local parametrization for $X_2$ (substituting in
\matrice):
$$\eqalign{
p_{02}=&{a_1^2-3a_0a_2\over a_0}v+{1\over2}w+\ldots\cr
p_{03}=&{\rm terms\ of\ degree}\ \ge2\cr
p_{12}=&-v+\ldots\cr
p_{13}=&{\rm terms\ of\ degree}\ \ge2}$$
\noindent This shows that $L$ represents a smooth point of $X_2$ and
its embedded tangent plane is $p_{03}=p_{13}=p_{23}=0$. i.e. the beta
plane $\beta(T_x\Sigma)$.
\qed

\proclaim Corollary \singflessi. If $d\ge4$, the congruence $X_2$
has a singular curve consisting of the closure of non-parabolic
inflectional lines meeting $\Sigma$ with multiplicity at least four.
The degree of this curve in $G(1,3)$ is $2d(d-3)(3d-2)$.

{\it Proof:} The first statement is an immediate corollary of Lemma
\tangenteflessi. The degree of the curve can be found in \Salmon,
art. 597 (page 286). An alternative way of computing this degree is to
use the construction in the proof of \bitangenti. The universal
quadratic form can be also viewed as a map $q^*Q^*(-1)\to Q$, so
that its determinant (whose zeros correspond to the pairs of
coincident points) is a section of $(\bigwedge^2Q)^{\otimes2}(2)$.
Intersecting $\tilde X_1$ with that class and the class of a
hyperplane one gets the wanted number. Of course, a better way would
be to work directly on ${\Bbb P}(Symm^4Q^*)$.
\qed

\bigskip

\noindent{\bf Remarks:} 1) From the invariants of the congruence of
bitangents $X_1$ found in Props. \invarianti\ and \bitangenti, the
degree of the (total) focal surface $F$ of $X_1$ must be
$d(d-3)(2d^3+2d^2-35d+26)$. Clearly, the strict focal surface $F_0$
is $\Sigma$. As already noticed in the proof of Prop. \classi, the
map $Y_1\to\Sigma$ has degree $(d+2)(d-3)$, so that $F_0$ counts with
multiplicity $(d+2)(d-3)$ in $F$. Therefore, $F$ has still some
extra components of total degree $2d(d-3)(d^3+d^2-18d+12)$.

2) Similarly, from the invariants of the congruence $X_2$ of flexes
to $\Sigma$ found in Prop. \invarianti, the degree of the total focal
surface $F$ of $X_2$ is $2d(6d^2-21d+16)$. The strict total surface
is again $\Sigma$. Since through a general point of $\Sigma$
there are two asymptotic lines and the ramification at each of them
is two (see Lemma \tangenteflessi), $\Sigma$ now counts with
multiplicity four. Hence, the extra components of $F$ have total
degree $2d(6d^2-21d+14)$.

The following propositions will explain where these extra components
come from.

\bigskip

\proclaim Proposition \rigatabitangenti. Let $\Sigma$ be a general
surface $\Sigma\subset{\Bbb P}^3$ of degree $d$ and let $X_1\subset
G(1,3)$ be the congruence of bitangents to $\Sigma$. Then there are
two curves of $X_1$ all of whose lines are entirely contained in the
(total) focal surface $F$ of $X_1$: The singular curve of Corollary
\singbitangenti\ and the curve of stationary bitangents to $\Sigma$
(i.e. bitangents such that the tangent plane to $\Sigma$ at the two
tangency points is the same). Moreover, the degree of the ruled
surface consisting of such stationary bitangents has degree
$d(d-2)(d-3)(d^2+2d-4)$.

{\it Proof:} Let $L$ be a bitangent tangent to $\Sigma$. If there is
only one tangency point, by Lemma \tangentebitangenti then $L$ has
intersection multiplicity at least four at the contact point and, if
this multiplicity is exactly four, then $L$ not contained in the focal
locus. But, if $\Sigma$ is general, the set of lines with intersection
multiplicity at least five at some point of $\Sigma$ should be finite
(and there would be precisely $5d(d-4)(7d-12)$ such lines). Hence there
is no curve of focal lines whose general element is tangent at two
infinitely close points.

Assume now that that there are two different tangency points
$x_1,\ x_2$. Suppose first that $L$ has intersection multiplicity at
least three at some of the points. Then, by Lemma
\tangentebitangenti, $L$ is a cuspidal point of $X_1$. Therefore, for
any point $x\in L$, the line counts at least twice as a line of the
congruence passing through $x$, which means that $L$ is entirely
contained in the focal locus.

So we assume that $L$ is simply tangent at $x_1$ and $x_2$, and let
$\Pi_1,\ \Pi_2$ be the respective embedded tangent planes to $\Sigma$.
Obviously the line $L$ (and hence also the points
$x_1$ and $x_2$) is contained in both $\Pi_1$ and $\Pi_2$. Then, by
Lemma \tangentebitangenti, the tangent plane to $X_1$ (as a surface in
${\Bbb P}^5$) at the point represented by $L$ is generated by the
pencils $\Omega(x_1,\Pi_2)$ and $\Omega(x_2,\Pi_1)$. Therefore, it is
clear that this plane is contained in $G(1,3)$ (and is in fact a
beta-plane) if and only if $\Pi_1=\Pi_2$.

Finally, the degree of the ruled surface of stationary bitangents
can be found in \Salmon, art. 613 (page 305). \qed

\bigskip

\noindent {\bf Remark:} Observe that $X_1$ possesses another
singular curve, namely the curve of tritangent lines. This is a
curve of degree ${1\over3}d(d-3)(d-4)(d-5)(d^2+3d-2)$, from
\Salmon, art. 599, pages 287-288. However, it is a triple nodal
curve (while the curve of Corollary \singbitangenti\ is a cuspidal
curve). This is what makes that its lines are not properly focal
lines.

\bigskip

\proclaim Proposition \rigataflessi. Let $\Sigma$ be a general
surface $\Sigma\subset{\Bbb P}^3$ of degree $d$ and let $X_2$ be the
congruences of flexes to $\Sigma$. Then there are two curves of $X_2$
all of whose lines are entirely contained in the (total) focal
surface $F$ of $X_2$: The singular curve of Corollary \singflessi\
and the curve of parabolic inflectional lines to $\Sigma$. Moreover,
the degree of the ruled surface of parabolic inflectional lines to
$\Sigma$ has degree $2d(d-2)(3d-4)$.

{\it Proof:} By Lemma \tangenteflessi, a curve consisting of focal
lines such that its general element is non-parabolic must be the
singular curve of asymptotic lines with intersection multiplicity at
least four. As in the previous Proposition \rigatabitangenti, that
curve clearly consists of focal lines.

Assume now that a general line of such a curve is parabolic. By Lemma
\tangenteflessi, a general point of such a curve (i.e. a line having
intersection multiplicity three at the tangency point) is a focal
line. The degree of the ruled surface of asymptotic lines at
parabolic points can be obtained as follows (of course, it can also
be found in \Salmon, art. 576, Ex. 3, page 255):

We observe from Lemma \tangentecomplesso\ that the surface $Y_2$
is double along its intersection with the surface $Y'$ (of pairs
$(x,L)\in Y$ with $x$ parabolic). Therefore, the degree of their
set-theoretical intersection will be ${1\over2}[Y_2][Y']\ell$. From
Prop. \classi, an easy calculation shows that the wanted degree is
$2d(d-4)(3d-4)$.
\qed

\bigskip
\bigskip

\noindent{\semilarge \S5. Smooth congruences of bitangents to
arbitrary surfaces in ${\Bbb P}^3$.}

\bigskip

In the previous section we dealt with congruences of bitangents and
flexes to smooth surfaces and, with the only exception of the
bitangents to a smooth quartic surface, we always got singular
congruences. However, our scope is to find smooth congruences. On
the other hand, we have seen that all lines of a smooth congruence are
bitangent to their focal surface, which is in general very singular.
So it is natural to study congruences of bitangents to arbitrary
surfaces in ${\Bbb P}^3$, hoping to then understand any smooth
congruence. The main problem is then how to compute the invariants of
of such a congruence. The bidegree is not difficult to find. We will
give it when the singularitites of $\Sigma$ and $\Sigma^*$ are not too
bad:

\proclaim Lemma \bigrado. Let $\Sigma\subset{\Bbb P}^3$ be a surface
of degree $d$, class $d^*$, class of the hyperplane section $\mu_1$,
ordinary nodal curve of degree $\delta$, ordinary cuspidal curve of
degree
$\kappa$ and no other singular curves. Assume the same hypothesis for
the singular locus of the dual surface holds, and let $\delta^*$ be the
number of bitangent planes through a point and $\kappa^*$ the number of
inflectional planes through a point. Then the bidegree of the
congruence $X$ of bitangents to $\Sigma$ is $(a,b)$ with $a={1\over
2}(\mu_1^2-3\kappa^*)+4d^*-5\mu_1$ and
$b={1\over2}(\mu_1^2-3\kappa)+4d-5\mu_1$.

{\it Proof:} The class $b$ is the number of lines of $X$ in a general
plane of ${\Bbb P}^3$, i.e. the number of bitangents of a general
hyperplane section of $\Sigma$. This hyperplane section has degree
$d$, class $\mu_1$, $\delta$ nodes and $\kappa$ cusps. Then, from
Pl\"ucker formulas (see for instance \Walker, V\S8.2) we get that
$d=\mu_1(\mu_1-1)-2b-3i$, $\kappa=3\mu_1(\mu_1-2)-6b-8i$ (where $i$
is the number of flexes of the curve). From this we immediately get
the wanted value for $b$. The value of $a$ is obtained by duality.
\qed

\bigskip

We start now a series of examples to try to illustrate what the
general situation should be.

\noindent{\bf Example \Veronese:} The hypothesis on the dual of
$\Sigma$ is really needed. For instance, consider the tangent
developable of a twisted cubic $C$. This is a a quartic surface
$\Sigma$ whose singular locus is $C$, which appears as a cuspidal
locus. Hence, $d=4$, $\mu_1=3$ (its hyperplane section is a rational
quartic with three cusps,  so its dual is a nodal cubic), $\delta=0$ and
$\kappa=3$. Then we get $b=1$ (in fact, as we remarked, the dual of the
hyperplane section of $\Sigma$ has one node).
But $\kappa^*=\delta^*=d^*=0$, since $\Sigma^*$ is a curve. Then the
formula for $a$ is not valid (fortunately, because the corresponding
value would be $a=-{21\over2}$, negative and not an integer!). The
correct value can be computed as follows.

Let $L$ be a bitangent line with tangency points $x_1$ and $x_2$. Then
obviously $L$ is the intersection of the tangent planes
$T_{x_1}\Sigma$ and $T_{x_2}\Sigma$. But the converse is also true.
Take two planes $\Pi_1$, $\Pi_2$ tangent to $\Sigma$. Since $\Sigma$ is
developable, they are tangent respectively along lines $L_1$, $L_2$.
Let $L$ be the intersection of $\Pi_1$ and $\Pi_2$. Then $L$ meets
$L_1$ in a point $x_1$ and meets $L_2$ in a point $x_2$. It is now
clear that $L$ is a bitangent line with contact points $x_1$ and
$x_2$. With this description, the dual congruence will be the
congruence of bisecants to the dual $\Sigma^*$ (which is a twisted
cubic). This dual congruence has bidegree $(1,3)$, so that our
congruence has bidegree $(3,1)$. Its total focal surface has degree
four (and a cuspidal curve), so it is precisely $\Sigma$ (contrary to
the situation for a smooth surface in ${\Bbb P}^3$, as we have seen in
Prop. \rigatabitangenti). Hence the congruence is the set of
bitangents to its focal surface (total or strict). This is not going
to be however the situation for a ``general'' congruence.

\bigskip

\noindent{\bf Example \seidue:} The above example shows that the
dual of the congruence of bisecants to a twisted cubic behaves
nicely with respect to is focal surface. So it is natural to see
what happens to the congruence $X$ dual of the other smooth
congruence of bisecants, namely the bisecants to an elliptic
quartic $C$. Then $X$ has bidegree $(a,b)=(6,2)$ and sectional
genus $g=3$. Hence, the total focal surface has degree $16$. On the
other hand, reasoning as in the previous example, $X$ will be the
congruence of bitangents to the dual $C^*$, which is a tangent
developable of degree $8$ and cuspidal curve of degree $12$
(corresponding to the osculating planes of $C^*$). Since the
hyperplane section has genus one, it follows easily that $\Sigma$
has a nodal curve of degree $\delta=8$ and hence $\mu_1=4$. What
happens now is that the total focal surface is twice $\Sigma=C^*$
(and therefore no formula for the invariants of the focal surface
is valid anymore). Indeed, given a general point $x\in\Sigma$,
there are two bitangents to $\Sigma$ with tangency points at $x$
and another point. Summing up, the congruence of bitangents to the
(strict) focal surface coincide with the congruence $X$ itself, but
the (total) focal surface of the congruence is not $\Sigma$ as a
scheme, but only as a set.

\bigskip

\noindent{\bf Example \dualWelters:} We have observed (Prop.
\bitangenti\ or Corollary \singbitangenti) that the only smooth
congruence of bitangents to a smooth surface in ${\Bbb P}^3$ is the
congruence of bidegree $(12,28)$ of bitangents to a smooth quartic
$\Sigma\subset{\Bbb P}^3$. By duality, we also have a smooth
congruence $X$ of bidegree $(28,12)$ consisting of the bitangents to
the dual $\Sigma^*$. This is a surface in ${\Bbb P}^3$ of degree $36$,
a nodal curve of degree $480$ and cuspidal curve of degree $96$. As in
the dual case, this counts six times in the total focal surface (since
through a general point of it there pass six lines that are tangent at
that point and another one). But the total focal surface has degree
$16$, so that there are no other components. Hence this congruence
verifies the same property with respect to the focal surface as the
one in the previous example.

\bigskip

\noindent{\bf Example \Kummer:} Consider in $G(1,3)$ the congruence
$X$ obtained in Example\duedue, which has bidegree $(2,2)$ and
sectional genus $g=1$. It is then a very classical result that the
focal surface is the so-called Kummer's surface, a quartic surface
with sixteen nodes, corresponding to the sixteen fundamental points
of $X$ (see for \Hudson\ for a thorough study of this surface).
However, the congruence of bisecants to the Kummer's surface (which
should have bidegree $(12,28)$) splits as sixteen beta-planes
(corresponding to the singular planes) and six congruences of
bidegree $(2,2)$ as above.

We conjecture that the general situation should be like the above
example (except for the existence of fundamental points). In other
words, a ``general congruence'' should have an irreducible reduced
focal surface (i.e. the total focal surface coincides with the strict
focal surface), and the congruence of bitangents to the focal surface
splits as the original congruence plus another congruence (in general
irreducible). Observe that the fact that the congruence of bitangents
to $F$ splits implies that one does not need to expect to have
excedentary components for the focal surface (as it should happen for
the congruence of all bitangents to a surface, as remarked in Prop.
\rigatabitangenti). Now we explicitly state our conjectures:

\proclaim Conjecture \focaleriducibile. If the total focal surface of a
smooth congruence $X$ is not irreducible, then either $X$ is the
congruence of secant lines to a curve in ${\Bbb P}^3$ (hence
necessarily the one in Example \duesei), or a congruence of bitangents
to a surface in ${\Bbb P}^3$ or a congruence of flexes to a surface in
${\Bbb P}^3$.

\proclaim Conjecture \focalenonridotta. It the total focal surface of
a smooth congruence $X$ is not reduced, the either $X$ is the
congruence of bitangents to a surface in ${\Bbb P}^3$ or a congruence
of flexes to a surface in ${\Bbb P}^3$.

These conjectures can be strengthen with the three following ones:

\proclaim Conjecture \bitangentiliscia. If the congruence of
bitangents to a surface $\Sigma\subset{\Bbb P}^3$ is smooth, then
either $\Sigma$ is a smooth quartic surface, or its dual (see
Example \dualWelters) or the tangent developable of a twisted cubic
(see Example \Veronese) or the one in Example \seidue.

\proclaim Conjecture \flessiliscia. There is no smooth congruence of
flexes to any surface in ${\Bbb P}^3$. More generally, there are no
congruences of the third class of Goldstein classification.

\proclaim Conjecture \bitangentifocale. Let $X$ be a smooth congruence
and let $F_0$ be its strict focal surface. Then $X$ coincides with the
congruence of bitangents to $F_0$ only in the case of Examples
\Veronese, \seidue, \dualWelters or its dual $(12,28)$ of bitangents
to a smooth quartic surface.
\bigskip
\bigskip

\noindent{\semilarge References:}
\bigskip

\item{\ArrondoGross} E. Arrondo -- M. Gross, {\it On smooth surfaces in
$Gr(1,{\bf P}^3)$ with a fundamental curve}, Manuscripta Math.,
79, (1993), 283-298.

\item{\Arkiv} E. Arrondo -- I. Sols -- R. Speiser, {\it Global moduli
of contacts}, Arkiv f\"or Math., 35 (1997), 1-57.

\item{\Ciro} C. Ciliberto -- E. Sernesi, {\it Singularities of the theta
divisor and congruences of planes}, Journal of Alg. Geom., 1 no. 2 (1992),
231-250.

\item{\Fano} G. Fano, {\it Studio di alcuni sistemi di rette
considerati come superficie dello spazio a cinque dimensioni}, Annali
di Matematica, 21 (1893), 141-192.

\item{\Goldstein} N. Goldstein, {\it The geometry of surfaces in
the 4-quadric}, Rend. Sem. Mat. Univers. Politecn. Torino, 43, 3
(1985), 467-499.

\item{\GH} P. Griffiths -- J. Harris, {\it Algebraic geometry and
local differential geometry}, Ann. Sci. \'Ecole Norm. Sup. (4) 12
(1979), 355-452.

\item{\Gross} M. Gross, {\it The distribution of bidegrees of smooth
surfaces in $G(1,{\Bbb P}^3)$}, Math. Ann. 292 (1992), 127-147.

\item{\Hudson} R. W. H. T. Hudson, {\it Kummer's quartic surface},
Cambridge Univ. Press, ed. 1990.

\item{\Johnsen} T. Johnsen, {\it Plane projections of a smooth
space curve}, in ``Parameter spaces'', Banach Center Publications,
VOl. 36 (1996), 89-110.

\item{\Schubert} S. Katz,  -- S.A. Str{\o}mme, {\it {\tt
schubert}, a {\tt Maple} package for intersection theory}, Available
at http://www.math.okstaste.edu/$\sim$katz/schubert.html or by
anonymous ftp from ftp.math.okstate.edu or linus.mi.uib.no, cd
pub/schubert.

\item{\McCrory} C. McCrory -- T. Shifrin, {\it Cusps of the projective
Gauss map}, J. Differential Geometry, 19 (1984), 257-276.

\item{\McCrorybis} C. McCrory -- T. Shifrin -- R. Varley, {\it The
Gauss map of a generic hypersurface in ${\Bbb P}^4$}, J.
Differential Geometry, 30 (1989), 689-759.

\item{\Peskine} C. Peskine -- L. Szpiro, {\it Liaison des
vari\'et\'es alg\'ebriques, I}, Invent. Math. 26 (1974), 271-302.

\item{\Rothuno} L. Roth, {\it Line congruences in three dimensions},
Proc. London Math. Soc. (2), 32 (1931), 72-86.

\item{\Rothdue} L. Roth, {\it Some properties of line congruences},
Proc. Camb. Phil. Soc., 27 (1931), 190-200.

\item{\Salmon} G. Salmon, {\it A treatise on the analytic geometry of
three dimension}, Vol. II, 5th ed. Chelsea Pub. Co., 1965.

\item{\Schumacher} R. Schumacher, {\it Classification der
algebraischen Strahlensysteme}, 37 (1890), 100-140.

\item{\Verra} A. Verra, {\it Geometria della retta in dimensione $2$}, 
unpublished paper (1986).

\item{\Walker} R. J. Walker, {\it Algebraic Curves}, Reprint by
Springer-Verlag, 1978.

\item{\Welters} G. E. Welters, {\it Abel-Jacobi isogenies for
certain types of Fano threefolds}, Mathematical Centre Tracts 141,
Amsterdam 1981.

\bigskip

\centerline{Authors address}

\bigskip

\centerline{Enrique Arrondo}
\centerline{Departamento de Algebra} \centerline{Facultad de
Ciencias Matem\'aticas} \centerline{Universidad Complutense de
Madrid} \centerline{28040 Madrid, Spain} \centerline{\tt
Enrique\_Arrondo@mat.ucm.es}

\bigskip

\centerline{Marina Bertolini and Cristina Turrini}
\centerline{Dipartimento di Matematica ``Federigo Enriques''}
\centerline{Universit\`a degli Studi di Milano} \centerline{Via C.
Saldini, 50} \centerline{20133 Milano, Italy} \centerline{\tt
Marina.Bertolini@mat.unimi.it \ \ \ \ \ \ \
Cristina.Turrini@mat.unimi.it}

\end